\documentclass[a4paper,12pt]{article}
\usepackage{amsmath}
\usepackage{latexsym}
\usepackage{url}
\usepackage{color}
\numberwithin{equation}{section}

\def\R{{\bf R}}

\def\N{{\bf N}}

\def\d{\displaystyle}
\def\e{{\varepsilon}}

\def\wt{\widetilde}

\def\p{\partial}
\def\v#1{\mbox{\boldmath $#1$}}

\newtheorem{thm}{Theorem}[section]

\newtheorem{prop}{Proposition}[section]
\newtheorem{rem}{Remark}[section]

\title{The generalized combined effect
 for one dimensional wave equations
with semilinear terms including product type}
\author{
Ryuki Kido
\footnote{
Master course, Mathematical Institute,
Tohoku University,
Aoba, Sendai 980-8578, Japan.
email: ryuki.kido.t1@dc.tohoku.ac.jp (Kido),
shu.takamatsu.r8@dc.tohoku.ac.jp (Takamatsu).},
Takiko Sasaki
\footnote{
Department of Mathematical Engineering, Faculty of Engineering, Musashino University,
3-3-3 Ariake, Koto-ku, Tokyo 135-8181, Japan./
Mathematical Institute, Tohoku University,
Aoba, Sendai 980-8578, Japan.
e-mail: t-sasaki@musashino-u.ac.jp.
},
Shu Takamatsu
\footnotemark[1]
,\\
Hiroyuki Takamura
\footnote{Mathematical Institute,
Tohoku University,
Aoba, Sendai 980-8578, Japan.
e-mail: hiroyuki.takamura.a1@tohoku.ac.jp.}
}
\date{{\small\it Dedicated to Professor Takayoshi Ogawa on his sixtieth birthday}
\[
\begin{array}{ll}
\mbox{\footnotesize{\bf Keywords:}}
& \mbox{\footnotesize semilinear wave equation, one dimension,
classical solution,}\\
& \mbox{\footnotesize lifespan, combined effect}\\
\mbox{\footnotesize{\bf MSC2020:}}
& \mbox{\footnotesize primary 35L71, secondary 35B44}\\
\end{array}
\]
}
\begin{document}
\maketitle
\begin{abstract}
We are interested in the so-called \lq\lq combined effect" of two different kinds
of nonlinear terms for semilinear wave equations in one space dimension.
Recently, the first result with the same formulation as in the higher dimensional case
has been obtained if and only if the total integral of the initial speed is zero,
namely Huygens' principle holds.
In this paper, we extend the nonlinear term to the general form including the product type. 
Such model equations are extremely meaningful only in one space dimension
because the most cases in higher dimensions possess the global-in-time existence
of a classical solution in the general theory for nonlinear wave equations.
It is also remarkable that our results on the lifespan estimates
are partially better than those of the general theory.
This fact tells us that there is a possibility to improve the general theory
which was expected complete more than 30 years ago.
\end{abstract}


\section{Introduction}

\par
Let us consider the initial value problems;
\begin{equation}
\label{IVP_gcombined}
\left\{
\begin{array}{ll}
	\d u_{tt}-u_{xx}=A|u_t|^p|u|^q+B|u|^r
	&\mbox{in}\quad \R\times(0,T),\\
	u(x,0)=\e f(x),\ u_t(x,0)=\e g(x),
	& x\in\R,
\end{array}
\right.
\end{equation}
where $p,q,r>1$, $A,B\ge0$ and $T>0$.
We assume that $f$ and $g$ are given smooth functions of compact support
and a parameter $\e>0$ is \lq\lq small enough".
We are interested in the lifespan $T(\e)$, the maximal existence time,
of classical solutions of (\ref{IVP_gcombined}).
Our results in this paper are the following estimates for $A>0$ and $B>0$;
\begin{equation}
\label{lifespan_non-zero}
T(\e)\sim
\min\{C\e^{-(p+q-1)},C\e^{-(r-1)/2}\}\qquad \mbox{if}\ \int_{\R}g(x)dx\not=0\\
\end{equation}
and
\begin{equation}
\label{lifespan_zero}
T(\e)\sim
\left\{
\begin{array}{l}
C\e^{-(p+q)(r-1)/(r+1)}\\
\quad\mbox{for}\ \d\frac{r+1}{2}\le p+q\le r,\\
\min\{C\e^{-(p+q-1)},C\e^{-r(r-1)/(r+1)}\}\\
\quad\mbox{otherwise}\\
\end{array}
\right.
 \mbox{if}\ \d\int_{\R}g(x)dx=0.
\end{equation}
Here we denote the fact that there are positive constants,
$C_1$ and $C_2$, independent of $\e$ satisfying $A(\e,C_1)\le T(\e)\le A(\e,C_2)$
by $T(\e)\sim A(\e,C)$.
We note that (\ref{lifespan_non-zero}) and (\ref{lifespan_zero})
are already established in the special setting $q=0$ by
Morisawa, Sasaki and Takamura \cite{MST},
but it is a non-trivial business to extend it to (\ref{IVP_gcombined})
due to the first term of product type for which different estimates from $q=0$
are required in the proof.
Also we note that the case of $p=1$, or $q=1$, is excluded
because there is no hope to construct a classical solution due to lack of the differentiability.
If we replace $|u_t|^p|u|^q$ with $u_t|u|^q$ for $p=1$ and $q>1$,
$|u_t|^pu$ for $p>1$ and $q=1$, $u_tu$ for $p=q=1$,
then we may have the similar result at least of the existence part,
but our method in this paper cannot be applicable directly for such terms.
\par
First we note that it was conjectured that
\[
T(\e)\sim C\e^{-(p+q-1)}\quad\mbox{for $A>0$ and $B=0$.}
\]
This was verified by Zhou \cite{Zhou} for the upper bound
with integer $p,q$ satisfying $p\ge1,q\ge0,p+q\ge2$,
and by Li,Yu and Zhou \cite{LYZ91,LYZ92}
for the lower bound with integer $p,q$ satisfying $p+q\ge2$ including more general
but smooth terms.
Note that \cite{Zhou} is a preprint version of Zhou \cite{Zhou01}
in which only the case of $q=0$ is considered.
But it is easy to apply its argument to the case of $q>0$ by making use of
\begin{equation}
\label{product}
|u_t|^p|u|^q=(p/(p+q))^p\left|\left(|u|^{(p+q)/p}\right)_t\right|^p
\end{equation}
as in \cite{Zhou}.
For the sake of completeness of this paper,
we shall repeat its proof in Appendix below.

\par
On the other hand, Zhou \cite{Zhou92} obtained
\[
T(\e)\sim
\left\{
\begin{array}{ll}
C\e^{-(r-1)/2} & \mbox{if}\ \d\int_{\R}g(x)dx\neq0,\\
C\e^{-r(r-1)/(r+1)} & \mbox{if}\ \d\int_{\R}g(x)dx=0
\end{array}
\right.
\mbox{for $A=0$ and $B>0$}.
\]
Therefore (\ref{lifespan_non-zero}) and (\ref{lifespan_zero}) are quite natural
as taking the minimum of both results
except for the first case in (\ref{lifespan_zero}),
in which, we have
\begin{equation}
\label{less_min}
\begin{array}{c}
\d C\e^{-(p+q)(r-1)/(r+1)}\le\min\{C\e^{-(p+q-1)},C\e^{-r(r-1)/(r+1)}\}\\
\d\quad\mbox{for}\ \d\frac{r+1}{2}\le p+q\le r.
\end{array}
\end{equation}
We shall call this special phenomenon by \lq\lq generalized combined effect" of two nonlinearities.
The original combined effect, which means the case of $q=0$, was first observed by
Han and Zhou \cite{HZ14} which targets to show the optimality of the result
of Katayama \cite{Katayama01} on the lower bound of the lifespan of classical solutions
of nonlinear wave equations with a nonlinear term $u_t^3+u^4$ in two space dimensions
including more general nonlinear terms.
It is known that $T(\e)\sim\exp\left(C\e^{-2}\right)$ for the nonlinear term
$u_t^3$ and $T(\e)=\infty$ for the nonlinear term $u^4$,
but Katayama \cite{Katayama01} obtained only a much worse estimate
than their minimum  as $T(\e)\ge c\e^{-18}$.
Surprisingly, more than ten years later, Han and Zhou \cite{HZ14} showed that
this result is optimal as $T(\e)\le C\e^{-18}$.
They also considered (\ref{IVP_gcombined}) with $q=0$ for all space dimensions $n$ bigger than 1
and obtain the upper bound of the lifespan.
Its counter part, the lower bound of the lifespan, was obtained by
Hidano, Wang and Yokoyama \cite{HWY16} for $n=2,3$.
See the introduction of \cite{HWY16} for the precise results and references.  
We note that the first case in (\ref{lifespan_zero}) with $q=0$ coincides with the lifespan estimate 
for the combined effect in \cite{HZ14, HWY16} if one sets $n=1$ formally.
Indeed, \cite{HZ14} and \cite{HWY16} showed that
\begin{equation}
\label{lifespan_highD}
T(\e)\sim C\e^{-2p(r-1)/\{2(r+1)-(n-1)p(r-1)\}}
\end{equation}
holds for $n=2,3$ provided
 \begin{equation}
 \label{condition_combined}
(r-1)\{(n-1)p-2\}<4,\ 2\le p\le r\le 2p-1,\ r>\frac{2}{n-1}.
\end{equation}
Later, Dai, Fang and Wang \cite{DFW19} improved the lower bound of lifespan
for the critical case in \cite{HWY16}.
They also show that $T(\e)<\infty$ for all $p,r>1$ in case of $n=1$, i.e. (\ref{IVP_gcombined})
with $q=0$.
For the non-Euclidean setting of the results above, see Liu and Wang \cite{LW20} for example,
in which the application to semilinear damped wave equations is included.

\par
Finally we strongly remark that our estimates in (\ref{lifespan_non-zero}) and (\ref{lifespan_zero})
are better than those of the general theory by Li, Yu and Zhou \cite{LYZ91, LYZ92}
in case of
\[
\frac{r+1}{2}<p+q<r\quad\mbox{and}\quad\int_{\R}g(x)dx=0
\]
with integer $p,q,r\ge2$.
Because our result on the lower bound of the lifespan can be established
also for the smooth terms as $u_{tt}-u_{xx}=u_t^pu^q+u^r$.
The typical example is $(p,q,r)=(2,2,6)$.
This fact shows a possibility to improve the general theory.
For details, see the last half of the next section.
We note that this kind of observations in Morisawa, Sasaki and Takamura \cite{MST}
has an error by wrong citation in the third case in (2.24) in \cite{MST}.
This paper corrects it.
We also note that, even for the original combined effect of $q=0$,
the integer points satisfying (\ref{condition_combined}) are
$(p,r)=(2,3),(3,3),(3,4)$ for $n=2$ and $(p,r)=(2,2)$ for $n=3$,
but (\ref{lifespan_highD}) with $p=r$ agrees with the case of $A=0$ and $B>0$.
See Introduction of Imai, Kato, Takamura and Wakasa \cite{IKTW}
for references on the case of $A=0$ and $B>0$.
Hence one can say that only the lifespan estimates with $(p,r)=(2,3),(3,4)$ for $n=2$
are essentially in the combined effect case.
If $q\neq0$, $p$ is replaced with $p+q$ in the results above.
Therefore it has less meaningful to consider (\ref{IVP_gcombined})
in higher space dimensions, $n\ge2$, if we discuss the optimality of the general theory. 
\par
Of course, some special structure of the nonlinear terms
such as \lq\lq null condition" guarantees the global-in-time existence.
See Nakamura \cite{Nakamura14}, Luli, Yang and Yu \cite{LYY18},
Zha \cite{Zha20, Zha22} for examples in this direction.
But we are interested in the optimality of the general theory.
The details are discussed at the end of Section 2 below.
This work is initiated by series of papers,
Kitamura \cite{Kitamura},
Kitamura, Morisawa and Takamura \cite{KMT22,KMT23},
Kitamura, Takamura and Wakasa \cite{KTW},
in which the weighted nonlinear terms are considered
for the purpose to be a trigger to extend the general theory
to the one for non-autonomous equations. 
\par
This paper is organized as follows.
In the next section, the preliminaries are introduced.
Moreover, (\ref{lifespan_non-zero}) and (\ref{lifespan_zero}) are divided into four theorems,
and we compare our results with those of the general theory.
Sections 3 is devoted to the proof of the existence part of (\ref{lifespan_non-zero}).
Sections 4 and 5 are devoted to the proof of the existence part of (\ref{lifespan_zero}).
Their main strategy is the iteration method in the weighted $L^\infty$ space 
due to Morisawa, Sasaki and Takamura \cite{MST}
which is originally introduced by John \cite{John79}.
Finally, we prove the blow-up part of (\ref{lifespan_non-zero}) and (\ref{lifespan_zero})
by following essentially, Han and Zhou \cite{HZ14} for the generalized combined effect,
and the iteration argument in \cite{MST} for other cases. 


\section{Preliminaries and main results}

Throughout this paper, we assume that the initial data
$(f,g)\in C_0^2(\R)\times C^1_0(\R)$ satisfies
\begin{equation}
\label{supp_initial}
\mbox{\rm supp }f,\ \mbox{supp }g\subset\{x\in\R:|x|\le R\},\quad R\ge1.
\end{equation}
Let $u$ be a classical solution of (\ref{IVP_gcombined}) in the time interval $[0,T]$.
Then the support condition of the initial data, (\ref{supp_initial}), implies that
\begin{equation}
\label{support_sol}
\mbox{supp}\ u(x,t)\subset\{(x,t)\in\R\times[0,T]:|x|\le t+R\}.
\end{equation}
For example, see Appendix of John \cite{John_book} for this fact.

\par
It is well-known that $u$ satisfies the following integral equation.
\begin{equation}
\label{u}
u(x,t)=\e u^0(x,t)+L(A|u_t|^p|u|^q+B|u|^r)(x,t),
\end{equation}
where $u^0$ is a solution of the free wave equation with the same initial data,
\begin{equation}
\label{u^0}
u^0(x,t):=\frac{1}{2}\{f(x+t)+f(x-t)\}+\frac{1}{2}\int_{x-t}^{x+t}g(y)dy,
\end{equation}
and a linear integral operator $L$ for a function $v=v(x,t)$ in Duhamel's term is defined by
\begin{equation}
\label{nonlinear}
L(v)(x,t):=\frac{1}{2}\int_0^tds\int_{x-t+s}^{x+t-s}v(y,s)dy.
\end{equation}
Then, one can apply the time-derivative to (\ref{u}) to obtain
\begin{equation}
\label{u_t}
u_t(x,t)=\e u_t^0(x,t)+L'(A|u_t|^p|u|^q+B|u|^r)(x,t)
\end{equation}
and
\begin{equation}
\label{u^0_t}
u_t^0(x,t)=\frac{1}{2}\{f'(x+t)-f'(x-t)+g(x+t)+g(x-t)\},
\end{equation}
where $L'$ for a function $v=v(x,t)$ is defined by
\begin{equation}
\label{nonlinear_derivative}
L'(v)(x,t):=\frac{1}{2}\int_0^t\{v(x+t-s,s)+v(x-t+s,s)\}ds.
\end{equation}
On the other hand, applying the space-derivative to (\ref{u}),
we have
\begin{equation}
\label{u_x}
u_x(x,t)=\e u_x^0(x,t)+\overline{L'}(A|u_t|^p|u|^q+B|u|^r)(x,t)
\end{equation}
and
\begin{equation}
\label{u^0_x}
u_x^0(x,t)=\frac{1}{2}\{f'(x+t)+f'(x-t)+g(x+t)-g(x-t)\},
\end{equation}
where $\overline{L'}$ for a function $v=v(x,t)$ is defined by
\begin{equation}
\label{nonlinear_derivative_conjugate}
\overline{L'}(v)(x,t):=
\frac{1}{2}\int_0^t\{v(x+t-s,s)-v(x-t+s,s)\}ds.
\end{equation}
Therefore, $u_x$ is expressed by $u$ and $u_t$.
Moreover, one more space-derivative to (\ref{u_t}) yields that
\begin{equation}
\label{u_tx}
\begin{array}{ll}
u_{tx}(x,t)=&\d\e u_{tx}^0(x,t)\\
&\d+AL'(p|u_t|^{p-2}u_tu_{tx}|u|^q+q|u|^{q-2}uu_x|u_t|^p)(x,t)\\
&\d+BL'(r|u|^{r-2}uu_x)(x,t)
\end{array}
\end{equation}
and
\begin{equation}
\label{u^0_tx}
u_{tx}^0(x,t)=\frac{1}{2}\{f''(x+t)-f''(x-t)+g'(x+t)+g'(x-t)\}.
\end{equation}
Similarly, we have that
\[
\begin{array}{ll}
u_{tt}(x,t)=&
\d\e u_{tt}^0(x,t)+A|u_t(x,t)|^p|u(x,t)|^q+B|u(x,t)|^r\\
&\d +A\overline{L'}(p|u_t|^{p-2}u_tu_{tx}|u|^q+q|u|^{q-2}uu_x|u_t|^p)(x,t)\\
&\d +B\overline{L'}(r|u|^{r-2}uu_x)(x,t)
\end{array}
\]
and
\[
u_{tt}^0(x,t)=\frac{1}{2}\{f''(x+t)+f''(x-t)+g'(x+t)-g'(x-t)\}.
\]
Therefore, $u_{tt}$ is expressed by $u,u_t,u_x,u_{tx}$ and so is $u_{xx}$, because of
\[
\begin{array}{ll}
u_{xx}(x,t)=
&\d\e u_{xx}^0(x,t)\\
&\d +A\overline{L'}(p|u_t|^{p-2}u_tu_{tx}|u|^q+q|u|^{q-2}uu_x|u_t|^p)(x,t)\\
&\d +B\overline{L'}(r|u|^{r-2}uu_x)(x,t)
\end{array}
\]
and
\[
u_{xx}^0(x,t)=u^0_{tt}(x,t).
\]

\par
First, we note the following fact.

\begin{prop}
\label{prop:system}
Assume that $(f,g)\in C^2(\R)\times C^1(\R)$.
Let $(u,w)$ be a $C^1$ solution of a system of integral equations;
\begin{equation}
\label{system}
\left\{
\begin{array}{l}
u=\e u^0+L(A|w|^p|u|^q+B|u|^r),\\
w=\e u_t^0+L'(A|w|^p|u|^q+B|u|^r)
\end{array}
\right.
\mbox{in}\ \R\times[0,T]
\end{equation}
with some $T>0$.
Then, $w\equiv u_t$ in $\R\times[0,T]$ holds
and $u$ is a classical solution of (\ref{IVP_gcombined}) in $\R\times[0,T]$.
\end{prop}
\par\noindent
{\bf Proof.} It is trivial that $w\equiv u_t$ by differentiating the first equation with respect to $t$.
 The rest part is easy along with the computations above in this section. 
\hfill$\Box$

\vskip10pt
Our results in (\ref{lifespan_non-zero}) and (\ref{lifespan_zero})
are divided into the following four theorems.

\begin{thm}
\label{thm:lower-bound_non-zero}
Let $A>0$ and $B>0$.
Assume (\ref{supp_initial}) and
\begin{equation}
\label{non-zero}
\int_{\R}g(x)dx\neq0.
\end{equation}
Then, there exists a positive constant $\e_1=\e_1(f,g,p,q,A,B,R)>0$ such that
a classical solution $u\in C^2(\R\times[0,T])$ of (\ref{IVP_gcombined}) exists
as far as $T$ satisfies
\begin{equation}
\label{lower-bound_non-zero}
T\le
\left\{
\begin{array}{ll}
c\e^{-(p+q-1)} & \mbox{for}\ p+q\le\d\frac{r+1}{2},\\
c\e^{-(r-1)/2} & \mbox{for}\ \d\frac{r+1}{2}\le p+q,
\end{array}
\right.
\end{equation}
where $0<\e\le\e_1$, and $c$ is a positive constant independent of $\e$.
\end{thm}

\begin{thm}
\label{thm:lower-bound_zero}
Let $A>0$ and $B>0$.
Assume (\ref{supp_initial}) and
\begin{equation}
\label{zero}
\int_{\R}g(x)dx=0.
\end{equation}
Then, there exists a positive constant $\e_2=\e_2(f,g,p,q,A,B,R)>0$ such that
a classical solution $u\in C^2(\R\times[0,T])$ of (\ref{IVP_gcombined}) exists
as far as $T$ satisfies
\begin{equation}
\label{lower-bound_zero}
T\le
\left\{
\begin{array}{ll}
c\e^{-(p+q-1)} & \mbox{for}\ p+q\le\dfrac{r+1}{2},\\
c\e^{-(p+q)(r-1)/(r+1)} & \mbox{for}\ \d\frac{r+1}{2}\le p+q\le r,\\
c\e^{-r(r-1)/(r+1)} & \mbox{for}\ p+q\ge r,
\end{array}
\right.
\end{equation}
where $0<\e\le\e_2$, and $c$ is a positive constant independent of $\e$.
\end{thm}

\begin{thm}
\label{thm:upper-bound_non-zero}
Let $A>0$ and $B>0$.
Assume (\ref{supp_initial}) and
\begin{equation}
\label{positive_non-zero}
\int_{\R}g(x)dx>0.
\end{equation}
Then, there exists a positive constant $\e_3=\e_3(f,g,p,q,A,B,R)>0$ such that
any classical solution of (\ref{IVP_gcombined}) in the time interval $[0,T]$ cannot exist
as far as $T$ satisfies
\begin{equation}
\label{upper-bound_non-zero}
T\ge
\left\{
\begin{array}{ll}
C\e^{-(p+q-1)} & \mbox{for}\ p+q\le\d\frac{r+1}{2},\\
C\e^{-(r-1)/2} & \mbox{for}\ \d\frac{r+1}{2}\le p+q,
\end{array}
\right.
\end{equation}
where $0<\e\le\e_3$, and $C$ is a positive constant independent of $\e$.
\end{thm}

\begin{thm}
\label{thm:upper-bound_zero}
Let $A>0$ and $B>0$.
Assume (\ref{supp_initial}) and
\begin{equation}
\label{positive_zero}
\begin{array}{l}
f(x)\ge0(\not\equiv0),\ g(x)\equiv0,\\
f(x)\ge f_0\ \mbox{and}\ -f'(x)\ge f_0\ \mbox{for}\ x\in(-R/2,0),
\end{array}
\end{equation}
with some positive constant $f_0$.
Then, there exists a positive constant $\e_4=\e_4(f,p,q,A,B,R)>0$ such that
any classical solution of (\ref{IVP_gcombined}) in the time interval $[0,T]$ cannot exist
as far as $T$ satisfies
\begin{equation}
\label{upper-bound_zero}
T\ge
\left\{
\begin{array}{ll}
C\e^{-(p+q-1)} & \mbox{for}\ p+q\le\dfrac{r+1}{2},\\
C\e^{-(p+q)(r-1)/(r+1)} & \mbox{for}\ \d\frac{r+1}{2}\le p+q\le r,\\
C\e^{-r(r-1)/(r+1)} & \mbox{for}\ p+q\ge r,
\end{array}
\right.
\end{equation}
where $0<\e\le\e_4$, and $C$ is a positive constant independent of $\e$.
\end{thm}

\begin{rem}
\label{rem:relation}
It is trivial that Theorem \ref{thm:lower-bound_non-zero} and
Theorem \ref{thm:upper-bound_non-zero} imply (\ref{lifespan_non-zero}).
On the other hand, we have that
\[
p+q-1=\frac{r(r-1)}{r+1}\quad\Longleftrightarrow\quad p+q=\frac{r^2+1}{r+1}
\]
and
\[
\frac{r+1}{2}<\frac{r^2+1}{r+1}<r.
\]
Moreover, we see that
\[
p+q-1\le\frac{(p+q)(r-1)}{r+1}\quad\Longleftrightarrow\quad p+q\le\frac{r+1}{2}.
\]
Therefore Theorem \ref{thm:lower-bound_zero} and
Theorem \ref{thm:upper-bound_zero} imply (\ref{lifespan_zero}).
\end{rem}

\par
The proofs of four theorems above appear in the following sections.
Form now on, we shall compare our results with those of the general theory 
by Li, Yu and Zhou \cite{LYZ91, LYZ92},
in which the following problem of general form is considered:
\begin{equation}
\label{IVP_general}
\left\{
\begin{array}{ll}
	\d u_{tt}-u_{xx}=F(u,Du,\p_xDu)
	&\mbox{in}\quad \R\times(0,\infty),\\
	u(x,0)=\e f(x),\ u_t(x,0)=\e g(x),
	& x\in\R,
\end{array}
\right.
\end{equation}
where we denote $D:=(\p_t,\p_x)$ and $F\in C^\infty(\R^5)$ satisfies
\[
F(\lambda)=O(|\lambda|^{1+\alpha})\quad\mbox{with $\alpha\in\N$ near $\lambda=0$}.
\]
(\ref{IVP_general}) requires $f,g\in C_0^\infty(\R)$.
Then, the lifespan of the classical solution of (\ref{IVP_general}) defined by $\wt{T}(\e)$ has
estimates from below as
\begin{equation}
\label{lifespan_general}
\wt{T}(\e)\ge
\left\{
\begin{array}{ll}
c\e^{-\alpha/2} & \mbox{in general},\\
c\e^{-\alpha(1+\alpha)/(2+\alpha)} & \mbox{if}\ \d\int_{\R}g(x)dx=0,\\
c\e^{-\min\{\beta_0/2,\alpha\}} & \mbox{if $\p_u^\beta F(0)=0$
for $1+\alpha\le\forall\beta\le\beta_0$}.
\end{array}
\right.
\end{equation}
This is the result of the general theory.
If one applies it to our problem (\ref{IVP_gcombined}) with
\begin{equation}
\label{F_special}
F(u,Du,\p_xDu)=u_t^pu^q+u^r\quad\mbox{with}\  p,q,r\in\N,
\end{equation}
one has the following estimates in each cases.
\begin{itemize}
\item
When $p+q<r$,
\par
then, we have to set $\alpha=p+q-1$ and $\beta_0=r-1$ which yield that
\[
\wt{T}(\e)\ge
\left\{
\begin{array}{ll}
c\e^{-(p+q-1)/2} & \mbox{in general},\\
c\e^{-(p+q)(p+q-1)/(p+q+1)} & \mbox{if}\ \d\int_{\R}g(x)dx=0,\\
c\e^{-\min\{(r-1)/2,p+q-1\}} &
\begin{array}{l}
\mbox{if $\p_u^\beta F(0)=0$}\\
\mbox{for $p+q\le\forall\beta\le r-1$}.
\end{array}
\end{array}
\right.
\]
We note that the third case is available for (\ref{F_special}).
Therefore, for $p+q\le(r+1)/2$, we obtain that
\[
\wt{T}(\e)\ge c\e^{-(p+q-1)}
\] 
whatever the value of $\d\int_{\R}g(x)dx$ is.
On the other hand, for $(r+1)/2<p+q$, i.e. 
\[
\frac{r-1}{2}<p+q-1,
\]
we obtain
\[
\wt{T}(\e)\ge
\left\{
\begin{array}{ll}
c\e^{-(r-1)/2} & \mbox{if}\ \d\int_{\R}g(x)dx\neq0,\\
c\e^{-\min\{(r-1)/2,(p+q)(p+q-1)/(p+q+1)\}} & \mbox{if}\ \d\int_{\R}g(x)dx=0.\\
\end{array}
\right.
\]
\item
When $p+q\ge r$,
\par
then, similarly to the case above, we have to set $\alpha=r-1$, which yields that
\[
\wt{T}(\e)\ge
\left\{
\begin{array}{ll}
c\e^{-(r-1)/2} & \mbox{in general},\\
c\e^{-r(r-1)/(r+1)} & \mbox{if}\ \d\int_{\R}g(x)dx=0,\\
c\e^{-\min\{\beta_0/2,(r-1)\}} & \mbox{if $\p_u^\beta F(0)=0$ for $r\le\forall\beta\le\beta_0$}.
\end{array}
\right.
\]
We note that the third case does not hold for (\ref{F_special}) by $\p_u^r F(0)\neq0$.
\end{itemize}
In conclusion, for the special nonlinear term in (\ref{F_special}),
the result of the general theory is
\[
\wt{T}(\e)\ge
\left\{
\begin{array}{ll}
c\e^{-(p+q-1)} & \mbox{for}\ p+q\le\d\frac{r+1}{2},\\
c\e^{-(r-1)/2} & \mbox{for}\ \d\frac{r+1}{2}\le p+q
\end{array}
\right.
\mbox{if}\ \int_{\R}g(x)dx\not=0\\
\]
and
\[
\begin{array}{c}
\wt{T}(\e)\ge
\left\{
\begin{array}{ll}
c\e^{-(p+q-1)} & \mbox{for}\ p+q\le\d\frac{r+1}{2},\\
c\e^{-\max\{(r-1)/2,(p+q)(p+q-1)/(p+q+1)\}} & \mbox{for}\ \d\frac{r+1}{2}\le p+q\le r,\\
c\e^{-r(r-1)/(r+1)} & \mbox{for}\ r\le p+q\\
\end{array}
\right.
\\
 \mbox{if}\ \d\int_{\R}g(x)dx=0.
 \end{array}
\]
Therefore a part of our results,
\begin{equation}
\label{better}
\begin{array}{l}
T(\e)\sim
C\e^{-(p+q)(r-1)/(r+1)}\\
\quad\mbox{if}\ \d\int_{\R}g(x)dx=0
\ \mbox{and}\
\frac{r+1}{2}< p+q<r,
\end{array}
\end{equation}
is better than the lower bound of $\wt{T}(\e)$ because of
\[
\frac{p+q-1}{p+q+1}<\frac{r-1}{r+1}.
\]
If one follows the proof in the following sections, one can find that it is easy to see that
our results on the lower bounds also hold for a special term (\ref{F_special})
by estimating the difference of nonlinear terms from above
after employing the mean value theorem.
We note that we have infinitely many examples of $(p,q,r)=(m,m,2m+1)$ as the inequality
\[
\frac{r+1}{2}=m+1<p+q=2m<r=2m+1
\]
holds for $m=2,3,4,\ldots$.
This fact indicates that we still have a possibility to improve the general theory
in the sense that the optimal results in (\ref{better}) should be included at least.


\section{Proof of Theorem \ref{thm:lower-bound_non-zero}}
\par
We basically employ the argument in Morisawa, Sasaki and Takamura \cite{MST} here.
According to Proposition \ref{prop:system},
we shall construct a $C^1$ solution of (\ref{system}).
Let $\{(u_j,w_j)\}_{j\in\N}$ be a sequence of $\{C^1(\R\times[0,T])\}^2$ defined by
\begin{equation}
\label{u_j,w_j}
\left\{
\begin{array}{ll}
u_{j+1}=\e u^0+L(A|w_j|^p|u_j|^q+B|u_j|^r), & u_1=\e u^0,\\
w_{j+1}=\e u_t^0+L'(A|w_j|^p|u_j|^q+B|u_j|^r), &  w_1=\e u_t^0.
\end{array}
\right.
\end{equation}
Then, in view of (\ref{u_x}) and (\ref{u_tx}), $\left((u_j)_x,(w_j)_x\right)$ has to satisfy
\begin{equation}
\label{u_j_x,w_j_x}
\left\{
\begin{array}{ll}
(u_{j+1})_x&=\e u^0_x+\overline{L'}\left(A|w_j|^p|u_j|^q+B|u_j|^r\right),\\
(u_1)_x&=\e u^0_x,\\
(w_{j+1})_x&=\e u^0_{tx}+L'\left(Ap|w_j|^{p-2}w_j(w_j)_x|u_j|^q\right)\\
&\quad+L'\left(Aq|u_j|^{q-2}u_j(u_j)_x|w_j|^p\right)\\
&\quad+L'\left(Br|u_j|^{r-2}u_j(u_j)_x\right),\\
 (w_1)_x&=\e u^0_{tx},
\end{array}
\right.
\end{equation}
so that the function space in which $\{(u_j,w_j)\}$ converges is
\[
\begin{array}{ll}
X:=&\{(u,w)\in\{C^1(\R\times[0,T])\}^2\ :\ \|(u,w)\|_X<\infty,\\
&\quad\mbox{supp}\ (u,w)\subset\{(x,t)\in\R\times[0,T]\ :\ |x|\le t+R\}\},
\end{array}
\]
which is equipped with a norm
\[
\|(u,w)\|_X:=\|u\|_1+\|u_x\|_1+\|w\|_2+\|w_x\|_2,
\]
where
\[
\begin{array}{l}
\d\|u\|_1:=\sup_{(x,t)\in\R\times[0,T]}|u(x,t)|,\\
\d\|w\|_2:=\sup_{(x,t)\in\R\times[0,T]}|(t-|x|+2R)w(x,t)|.
\end{array}
\]
First we note that supp $(u_j,w_j)\subset\{(x,t)\in\R\times[0,T]\ :\ |x|\le t+R\}$ implies supp
$(u_{j+1},w_{j+1})\subset\{(x,t)\in\R\times[0,T]\ :\ |x|\le t+R\}$.
It is easy to check this fact by assumption on the initial data (\ref{supp_initial})
and the definitions of $L,\overline{L},L',\overline{L'}$ in the previous section.

\par
The following lemma contains some useful a priori estimates.
\begin{prop}
\label{prop:apriori}
Let $(u,w)\in\{C(\R\times[0,T])\}^2$ and supp\ $(u,w)\subset\{(x,t)\in\R\times[0,T]:|x|\le t+R\}$. Then there exists a positive constant $C$ independent of $T$ and $\e$ such that
\begin{equation}
\label{apriori}
\begin{array}{l}
\|L(|w|^p|u|^q)\|_1\le C\|w\|_2^p\|u\|_1^q(T+R),\ \|L(|u|^r)\|_1\le C\|u\|_1^r(T+R)^2,\\
\|L'(|w|^p|u|^q)\|_2\le C\|w\|_2^p\|u\|_1^q(T+R),\ \|L'(|u|^r)\|_2\le C\|u\|_1^r(T+R)^2,\\
\|L'(|w|^p|u|^q)\|_1\le C\|w\|_2^p\|u\|_1^q(T+R),\ \|L'(|u|^r)\|_1\le C\|u\|_1^r(T+R)^2.
\end{array}
\end{equation}
\end{prop}

\par\noindent
{\bf Proof.} The proof of Proposition \ref{prop:apriori} is completely same as
the one of Proposition 3.1 in Morisawa, Sasaki and Takamura \cite{MST} because
$\|u\|_1$ has no weight, so that $\|w\|_2^p$ in \cite{MST} is simply replaced with
$\|w\|_2^p\|u\|_1^q$.
\hfill$\Box$

\par
Let us continue to prove Theorem \ref{thm:lower-bound_non-zero}.
Set
\[
M:=\sum_{\alpha=0}^2\|f^{(\alpha)}\|_{L^\infty(\R)}
+\|g\|_{L^1(\R)}+\sum_{\beta=0}^1\|g^{(\beta)}\|_{L^\infty(\R)}.
\]
\vskip10pt
\par\noindent
{\bf The convergence of the sequence $\v{\{(u_j,w_j)\}}$.}
\par
First we note that $\|u_1\|_1,\|w_1\|_2\le M\e$ by (\ref{u^0}) and (\ref{u^0_t}).
Since (\ref{u_j,w_j}) and (\ref{apriori}) yield that
\[
\left\{
\begin{array}{ll}
\|u_{j+1}\|_1
&\le M\e+A\|L(|w_j|^p|u_j|^q)\|_1+B\|L(|u_j|^r))\|_1\\
&\le M\e+AC\|w_j\|_2^p\|u_j\|_1^q(T+R)+BC\|u_j\|^r_1(T+R)^2,\\
\|w_{j+1}\|_2
&\le M\e+A\|L'(|w_j|^p|u_j|^q)\|_2+B\|L'(|u_j|^r)\|_2\\
&\le M\e+AC\|w_j\|_2^p\|u_j\|_1^q(T+R)+BC\|u_j\|^r_1(T+R)^2,\
\end{array}
\right.
\]
the boundedness of $\{(u_j,w_j)\}$, i.e.
\begin{equation}
\label{bound_(u,w)}
\|u_j\|_1,\|w_j\|_2\le 3M\e\quad(j\in\N),
\end{equation}
follows from
\begin{equation}
\label{condi1}
AC(3M\e)^{p+q}(T+R),BC(3M\e)^r(T+R)^2\le M\e.
\end{equation}
Assuming (\ref{condi1}), one can estimate $(u_{j+1}-u_j)$ and $(w_{j+1}-w_j)$ as follows.
Making use of
\[
\begin{array}{l}
|w_j|^p|u_j|^q-|w_{j-1}|^p|u_{j-1}|^q\\
=(|w_j|^p-|w_{j-1}|^p)|u_j|^q+|w_{j-1}|^p(|u_j|^q-|u_{j-1}|^q)
\end{array}
\]
and
\[
\left||u_j|^r-|u_{j-1}|^r\right|\le2^{r-1}r\left(|u_j|^{r-1}+|u_{j-1}|^{r-1}\right)|u_j-u_{j-1}|
\]
together with Proposition \ref{prop:apriori}, we have that
\[
\begin{array}{l}
\|u_{j+1}-u_j\|_1\\
\le\|L(A|w_j|^p|u_j|^q-A|w_{j-1}|^p|u_{j-1}|^q+B|u_j|^r-B|u_{j-1}|^r)\|_1\\
\le 2^{p-1}pA\|L\left((|w_j|^{p-1}+|w_{j-1}|^{p-1})|w_j-w_{j-1}||u_j|^q\right)\|_1\\
\quad+2^{q-1}qA\|L\left(|w_{j-1}|^p(|u_j|^{q-1}+|u_{j-1}|^{q-1})|u_j-u_{j-1}|\right)\|_1\\
\quad+2^{r-1}rB\|L\left((|u_j|^{r-1}+|u_{j-1}|^{r-1})|u_j-u_{j-1}|\right)\|_1\\
\le 2^{p-1}pAC(T+R)(\|w_j\|_2^{p-1}+\|w_{j-1}\|_2^{p-1})\|u_j\|_1^q\|w_j-w_{j-1}\|_2\\
\quad+2^{q-1}qAC(T+R)\|w_{j-1}\|_2^p(\|u_j\|_1^{q-1}+\|u_{j-1}\|_1^{q-1})\|u_j-u_{j-1}\|_1\\
\quad+2^{r-1}rBC(T+R)^2(\|u_j\|_1^{r-1}+\|u_{j-1}\|_1^{r-1})\|u_j-u_{j-1}\|_1\\
\le 2^ppAC(3M\e)^{p+q-1}(T+R)\|w_j-w_{j-1}\|_2\\
\quad+2^qqAC(3M\e)^{p+q-1}(T+R)\|u_j-u_{j-1}\|_1\\
\quad+2^rrBC(3M\e)^{r-1}(T+R)^2\|u_j-u_{j-1}\|_1
\end{array}
\]
and similarly
\[
\begin{array}{ll}
\|w_{j+1}-w_j\|_2
&\le 2^ppAC(3M\e)^{p+q-1}(T+R)\|w_j-w_{j-1}\|_2\\
&\quad+2^qqAC(3M\e)^{p+q-1}(T+R)\|u_j-u_{j-1}\|_1\\
&\quad+2^rrBC(3M\e)^{r-1}(T+R)^2\|u_j-u_{j-1}\|_1.
\end{array}
\]
Here we employ H\"older's inequality to obtain
\[
\begin{array}{l}
\|L(|w_j|^{p-1}|w_j-w_{j-1}|)|u_j|^q\|_1\\
=\|L\left(\left||w_j|^{(p-1)/p}|w_j-w_{j-1}|^{1/p}\right|^p|u_j|^q\right)\|_1\\
\le C(T+R)\||w_j|^{(p-1)/p}|w_j-w_{j-1}|^{1/p}\|_2^p\|u_j\|_1^q\\
\le C(T+R)\|w_j\|_2^{p-1}\|w_j-w_{j-1}\|_2\|u_j\|_1^q
\end{array}
\]
and so on.
Therefore the convergence of $\{u_j\}$ follows from
\begin{equation}
\label{convergence}
\begin{array}{l}
\|u_{j+1}-u_j\|_1+\|w_{j+1}-w_j\|_2\\
\d\le\frac{1}{2}\left(\|u_j-u_{j-1}\|_1+\|w_j-w_{j-1}\|_2\right)
\end{array}
\quad(j\ge2)
\end{equation}
provided (\ref{condi1}) and
\begin{equation}
\label{condi2}
2^{p+q}(p+q)AC(3M\e)^{p+q-1}(T+R),2^rrBC(3M\e)^{r-1}(T+R)^2\le\frac{1}{8}
\end{equation}
are fulfilled.

\vskip10pt
\par\noindent
{\bf The convergence of the sequence $\v{\{\left((u_j)_x,(w_j)_x\right)\}}$.}
\par
First we note that $\|(u_1)_x\|_1,\|(w_1)_x\|_2\le M\e$ by (\ref{u^0_x}) and (\ref{u^0_tx}).
Assume that (\ref{condi1}) and (\ref{condi2}) are fulfilled.
Since (\ref{u_j_x,w_j_x}) and (\ref{apriori}) yield that
\[
\begin{array}{ll}
\|(u_{j+1})_x\|_1
&\le M\e+A\|\overline{L'}\left(|w_j|^p|u_j|^q|\right)\|_1+B\|\overline{L'}\left(|u_j|^r|\right)\|_1\\
&\le M\e+AC(T+R)\|w_j\|_2^p\|u_j\|_1^q+BC(T+R)^2\|u_j\|_1^r\\
&\le M\e+AC(3M\e)^{p+q}(T+R)+BC(3M\e)^r(T+R)^2
\end{array}
\]
because of a trivial property  $|\overline{L'}(v)|\le L'(|v|)$ and
\[
\begin{array}{ll}
\|(w_{j+1})_x\|_2
&\le M\e+pA\|L'\left(|w_j|^{p-1}|(w_j)_x||u_j|^q\right)\|_2\\
&\quad+qA\|L'\left(|w_j|^p|u_j|^{q-1}|(u_j)_x|\right)\|_2\\
&\quad+rB\|L'\left(|u_j|^{r-1}|(u_j)_x|\right)\|_2\\
&\le M\e+pAC(T+R)\|w_j\|_2^{p-1}\|(w_j)_x\|_2\|u_j\|_1^q\\
&\quad+qAC(T+R)\|w_j\|_2^p\|u_j\|_1^{q-1}\|(u_j)_x\|_1\\
&\quad+rBC(T+R)^2\|u_j\|_1^{r-1}\|(u_j)_x\|_1\\
&\le M\e+pAC(3M\e)^{p+q-1}(T+R)\|(w_j)_x\|_2\\
&\quad+qAC(3M\e)^{p+q-1}(T+R)\|(u_j)_x\|_1\\
&\quad+rBC(3M\e)^{r-1}(T+R)^2\|(u_j)_x\|_1.
\end{array}
\]
The boundedness of $\{\left((u_j)_x,(w_j)_x\right)\}$, i.e.
\begin{equation}
\label{bound_U_x}
\|(u_j)_x\|_1,\|(w_j)_x\|_2\le 3M\e\quad(j\in\N),
\end{equation}
follows from
\begin{equation}
\label{condi3}
(p+q)AC(3M\e)^{p+q}(T+R),rBC(3M\e)^r(T+R)^2\le M\e.
\end{equation}
Assuming (\ref{condi3}), one can estimate $\{(u_{j+1})_x-(u_j)_x\}$
and $\{(w_{j+1})_x-(w_j)_x\}$ as follows.
It is easy to see that
\[
\begin{array}{ll}
\|(u_{j+1})_x-(u_j)_x\|_1
&\le A\|\overline{L'}(|w_j|^p|u_j|^q-|w_{j-1}|^p|u_{j-1}|^q)\|_1\\
&\quad+B\|\overline{L'}(|u_j|^r-|u_{j-1}|^r)\|_1,
\end{array}
\]
which can be handled like $(w_{j+1}-w_j)$ as before, so that we have that
\[
\begin{array}{ll}
\|(u_{j+1})_x-(u_j)_x\|_1
&\le 2^ppAC(3M\e)^{p+q-1}(T+R)\|w_j-w_{j-1}\|_2\\
&\quad+2^qqAC(3M\e)^{p+q-1}(T+R)\|u_j-u_{j-1}\|_1\\
&\quad+2^rrBC(3M\e)^{r-1}(T+R)^2\|u_j-u_{j-1}\|_1
\end{array}
\]
because of $|\overline{L'}(v)|\le L'(|v|)$, which implies that
\begin{equation}
\label{convergence_u_x}
\|(u_{j+1})_x-(u_j)_x\|_1
=O\left(\frac{1}{2^j}\right)
\end{equation}
as $j\rightarrow\infty$ due to (\ref{convergence}).

\par
On the other hand, we have that
\[
\begin{array}{l}
\|(w_{j+1})_x-(w_j)_x\|_2\\
\le pA\|L'(|w_j|^{p-2}w_j(w_j)_x|u_j|^q-|w_{j-1}|^{p-2}w_{j-1}(w_{j-1})_x)|u_{j-1}|^q\|_2\\
\quad+qA\|L'(|w_j|^p|u_j|^{q-2}u_j(u_j)_x-|w_{j-1}|^p|u_{j-1}|^{q-2}u_{j-1}(u_{j-1})_x)\|_2\\
\quad+rB\|L'(|u_j|^{r-2}u_j(u_j)_x-|u_{j-1}|^{r-2}u_{j-1}(u_{j-1})_x)\|_2.
\end{array}
\]
The first term on the right hand side of this inequality is divided into three pieces
according to
\[
\begin{array}{l}
|w_j|^{p-2}w_j(w_j)_x|u_j|^q-|w_{j-1}|^{p-2}w_{j-1}(w_{j-1})_x|u_{j-1}|^q\\
=(|w_j|^{p-2}w_j-|w_{j-1}|^{p-2}w_{j-1})(w_j)_x|u_j|^q\\
\quad+|w_{j-1}|^{p-2}w_{j-1}((w_j)_x-(w_{j-1})_x)|u_j|^q\\
\quad+|w_{j-1}|^{p-2}w_{j-1}(w_{j-1})_x(|u_j|^q-|u_{j-1}|^q).
\end{array}
\]
Since one can employ the estimate
\[
\begin{array}{l}
\left||w_j|^{p-2}w_j-|w_{j-1}|^{p-2}w_{j-1}\right|
\\
\le
\left\{
\begin{array}{ll}
(p-1)2^{p-2}(|w_j|^{p-2}+|w_{j-1}|^{p-2})|w_j-w_{j-1}| & \mbox{for}\ p\ge2,\\
2|w_j-w_{j-1}|^{p-1} & \mbox{for}\ 1<p<2,
\end{array}
\right.
\end{array}
\]
and the same one in which $w$ is replaced with $u$, we obtain that
\[
\begin{array}{l}
\|(w_{j+1})_x-(w_j)_x\|_2\\
\le pAC(T+R)\|(w_j)_x\|_2\|u_j\|_1^q\times\\
\qquad\times
\left\{
\begin{array}{ll}
(p-1)2^{p-2}(\|w_j\|_2^{p-2}+\|w_{j-1}\|_2^{p-2})\|w_j-w_{j-1}\|_2 & \mbox{for}\ p\ge2,\\
2\|w_j-w_{j-1}\|_2^{p-1} & \mbox{for}\ 1<p<2
\end{array}
\right.
\\
\quad+pAC(T+R)\|w_{j-1}\|_2^{p-1}\|(w_j)_x-(w_{j-1})_x\|_2\|u_j\|_1^q\\
\quad+2^qpqAC(T+R)\|w_{j-1}\|_2^{p-1}\|(w_{j-1})_x\|_2\times\\
\qquad\times(\|u_j\|_1^{q-1}+\|u_{j-1}\|_1^{q-1})\|u_j-u_{j-1}\|_1\\
\quad+qAC(T+R)\|(u_j)_x\|_1\|w_j\|_2^p\times\\
\qquad\times
\left\{
\begin{array}{ll}
(q-1)2^{q-2}(\|u_j\|_1^{q-2}+\|u_{j-1}\|_1^{q-2})\|u_j-u_{j-1}\|_1 & \mbox{for}\ q\ge2,\\
2\|u_j-u_{j-1}\|_1^{q-1} & \mbox{for}\ 1<q<2
\end{array}
\right.
\\
\quad+qAC(T+R)\|u_{j-1}\|_1^{q-1}\|(u_j)_x-(u_{j-1})_x\|_1\|w_j\|_2^p\\
\quad+2^ppqAC(T+R)\|u_{j-1}\|_1^{q-1}\|(u_{j-1})_x\|_1\times\\
\qquad\times(\|w_j\|_2^{p-1}+\|w_{j-1}\|_2^{p-1})\|w_j-w_{j-1}\|_2\\
\quad +rBC(T+R)^2\|(u_j)_x\|_1\times\\
\qquad\times
\left\{
\begin{array}{ll}
(r-1)2^{r-2}(\|u_j\|_1^{r-2}+\|u_{j-1}\|_1^{r-2})\|u_j-u_{j-1}\|_1 & \mbox{for}\ r\ge2,\\
2\|u_j-u_{j-1}\|_1^{r-1} & \mbox{for}\ 1<r<2
\end{array}
\right.
\\
\quad+rBC(T+R)^2\|u_{j-1}\|_1^{r-1}\|(u_j)_x-(u_{j-1})_x\|_1.
\end{array}
\]
Hence it follows from (\ref{convergence}) and (\ref{convergence_u_x}) that
\[
\begin{array}{ll}
\|(w_{j+1})_x-(w_j)_x\|_2
&\le pAC(3M\e)^{p+q-1}(T+R)\|(w_j)_x-(w_{j-1})_x\|_2\\
&\quad\d+
O\left(\frac{1}{2^{j\min\{p-1,q-1,r-1,1\}}}\right) 
\end{array}
\]
as $j\rightarrow\infty$.
Therefore we obtain the convergence of $\{\left((u_j)_x,(w_j)_x\right)\}$ provided
\begin{equation}
\label{condi4}
pAC(3M\e)^{p+q-1}(T+R)\le\frac{1}{2}.
\end{equation}

\vskip10pt
\par\noindent
{\bf Continuation of the proof.}
\par
The convergence of the sequence $\{(u_j,w_j)\}$ to $(u,w)$ in the closed subspace of $X$
satisfying
\[
 \|u\|_1,\|(u_x)\|_1,\|w\|_2,\|(w)_x\|_2\le 3M\e
 \]
is established by (\ref{condi1}), (\ref{condi2}), (\ref{condi3}), and (\ref{condi4}),
which follow from
\[
C_0\e^{p+q-1}(T+R)\le1\quad\mbox{and}\quad C_0\e^{r-1}(T+R)^2\le1,
\]
where
\[
\begin{array}{ll}
C_0:=\max &
\{3^{p+q}ACM^{p+q-1},3^rBCM^{r-1},\\
&\quad2^{p+q+3}3^{p+q-1}(p+q)ACM^{p+q-1}, 2^{r+3}3^{r-1}rBCM^{r-1},\\
&\quad 3^{p+q}(p+q)ACM^{p+q-1},3^rrBCM^{r-1},\\
&\quad2\cdot3^{p+q-1}pACM^{p+q-1}\}.
\end{array}
\]
Therefore the statement of Theorem \ref{thm:lower-bound_non-zero}
is established with
\[
\left\{
\begin{array}{l}
\d c=\frac{1}{2}\min\{C_0^{-1},C_0^{-1/2}\},\\
\e_1:=\min\{(2C_0R)^{-1/(p+q-1)},(2^2C_0R^2)^{-1/(r-1)}\}
\end{array}
\right.
\]
because
\[
R\le\frac{1}{2}\min\{C_0^{-1}\e^{-(p+q-1)},C_0^{-1/2}\e^{-(r-1)/2}\}
\]
holds for $0<\e\le\e_1$.
\hfill$\Box$


\section{Proof of Theorem \ref{thm:lower-bound_zero}}
\par
In this section also, 
we basically employ the argument in Morisawa, Sasaki and Takamura \cite{MST}.
But the different estimates for the product term are required here. 
First we note that the strong Huygens' principle
\begin{equation}
\label{Huygens}
u^0(x,t)\equiv0\quad\mbox{in}\ D
\end{equation}
holds in this case of (\ref{zero}), where
\[
D:=\{(x,t)\in\R\times[0,\infty)\ :\ t-|x|\ge R\}.
\]
This is almost trivial if one takes a look on the representation of $u^0$ in (\ref{u^0})
and the support condition on the data in (\ref{supp_initial}).
But one can see it also by Proposition 2.2
in Kitamura, Morisawa and Takamura \cite{KMT22} for the details.
So, our unknown functions are
$U:=u-\e u^0$ and $W:=w-\e u_t^0$ in (\ref{system}).
Let $\{(U_j,W_j)\}_{j\in\N}$ be a sequence of $\{C^1(\R\times[0,T])\}^2$ defined by
\begin{equation}
\label{U_j,W_j}
\left\{
\begin{array}{ll}
U_{j+1}=L(A|W_j+\e u_t^0|^p|U_j+\e u^0|^q+B|U_j+\e u^0|^r), & U_1=0,\\
W_{j+1}=L'(A|W_j+\e u_t^0|^p|U_j+\e u^0|^q+B|U_j+\e u^0|^r), &  W_1=0.
\end{array}
\right.
\end{equation}
Then, $\left\{\left((U_j)_x,(W_j)_x\right)\right\}$ has to satisfy
\begin{equation}
\label{U_j_x,W_j_x}
\left\{
\begin{array}{ll}
(U_{j+1})_x&=\overline{L'}\left(A|W_j+\e u_t^0|^p|U_j+\e u^0|^q+B|U_j+\e u^0|^r\right),\\
(U_1)_x &=0,\\
(W_{j+1})_x 
& =L'\left(Ap|W_j+\e u_t^0|^{p-2}(W_j+\e u_t^0)((W_j)_x+\e u_{tx}^0)|U_j+\e u^0|^q\right)\\
& \quad+L'\left(Aq|W_j+\e u_t^0|^p|U_j+\e u^0|^{q-2}(U_j+\e u^0)((U_j)_x+\e u_x^0)\right)\\
& \quad+L'\left(Br|U_j+\e u^0|^{r-2}(U_j+\e u^0)((U_j)_x+\e u_x^0)\right),\\
(W_1)_x &=0,
\end{array}
\right.
\end{equation}
so that the function space in which $\{(U_j,W_j)\}$ converges is
\[
\begin{array}{ll}
Y:=&\{(U,W)\in\{C^1(\R\times[0,T])\}^2\ :\ \|(U,W)\|_Y<\infty,\\
&\quad \mbox{supp}\ (U,W)\subset\{|x|\le t+R\}\}
\end{array}
\]
which is equipped with a norm 
\[
\|(U,W)\|_Y:=\|U\|_3+\|U_x\|_3+\|W\|_4+\|W_x\|_4,
\]
where
\[
\begin{array}{l}
\d\|U\|_3:=\sup_{(x,t)\in\R\times[0,T]}(t+|x|+R)^{-1}|U(x,t)|,\\
\d\|W\|_4:=\sup_{(x,t)\in\R\times[0,T]}\{\chi_D(x,t)+(1-\chi_D(x,t))(t+|x|+R)^{-1}\}|W(x,t)|,
\end{array}
\]
and $\chi_D$ is a characteristic function of $D$.
 Similarly to the proof of Theorem \ref{thm:lower-bound_non-zero},
we note that supp $(U_j,W_j)\in\{(x,t)\in\R\times[0,T]\ :\ |x|\le t+R\}$ implies
supp $(U_{j+1},W_{j+1})\in\{(x,t)\in\R\times[0,T]\ :\ |x|\le t+R\} $.

\par
The following lemmas are a priori estimates in this case.
\begin{prop}
\label{prop:apriori_linear}
Let $(U,W)\in\{C(\R\times[0,T])\}^2$ with
\[
\mbox{\rm supp}\ (U,W)\subset\{(x,t)\in\R\times[0,T]:|x|\le t+R\}
\]
and $U^0\in C(\R\times[0,T])$ with
\[
\mbox{\em supp}\ U^0\subset\{(x,t)\in\R\times[0,T]:(t-|x|)_+\le|x|\le t+R\}.
\]
Then there exists a positive constant $E$ independent of $T$ such that
\begin{equation}
\label{apriori_linear}
\left\{
\begin{array}{ll}
\|L(|U^0|^{q-m}|W|^m)\|_3 & \le E\|U^0\|_\infty^{q-m}\|W\|_4^m(T+R)^m,\\
\|L(|U^0|^{p-m}|U|^m)\|_3 & \le E\|U^0\|_\infty^{p-m}\|U\|_3^m(T+R)^m,\\
\|L'(|U^0|^{q-m}|W|^m)\|_4 & \le E\|U^0\|_\infty^{q-m}\|W\|_4^m(T+R)^m,\\
\|L'(|U^0|^{p-m}|U|^m)\|_4 & \le E\|U^0\|_\infty^{p-m}\|U\|_3^m(T+R)^m,\\
\|L'(|U^0|^{q-m}|W|^m)\|_3 & \le E\|U^0\|_\infty^{q-m}\|W\|_4^m(T+R)^m,\\
\|L'(|U^0|^{p-m}|U|^m)\|_3 & \le E\|U^0\|_\infty^{p-m}\|U\|_3^m(T+R)^m,
\end{array}
\right.
\end{equation}
where $p-m,q-m>0\ (m=0,1,2)$ and the norm $\|\cdot\|_\infty$ is defined by
\[
\|U^0\|_\infty:=\sup_{(x,t)\in\R\times[0,T]}|U^0(x,t)|.
\]
\end{prop}

\par\noindent
{\bf Proof.} This lemma is exactly same as Proposition 5.1 in Morisawa, Sasaki and Takamura \cite{MST}.
\hfill$\Box$

\vskip10pt
\par
We note that $U^0$ in the theorem above will be replaced with
$u^0$ or $u_t^0$ and their spatial derivatives later
due to (\ref{Huygens}).

\begin{prop}
\label{prop:apriori_zero}
Let $(U,W)\in\{C(\R\times[0,T])\}^2$ with
\[
\mbox{\rm supp}\ (U,W)\subset\{(x,t)\in\R\times[0,T]:|x|\le t+R\}.
\]
Then there exists a positive constant $C$ independent of $T$ such that
\begin{equation}
\label{apriori_zero}
\left\{
\begin{array}{ll}
\|L(|W|^p|U|^q)\|_3 & \le C\|W\|_4^p\|U\|_3^q(T+R)^{p+q},\\
\|L(|U|^r)\|_3 & \le C\|U\|_3^r(T+R)^{r+1},\\
\|L'(|W|^p|U|^q)\|_4 & \le C\|W\|_4^p\|U\|_3^q(T+R)^{p+q},\\
\|L'(|U|^r)\|_4 & \le C\|U\|_3^r(T+R)^{r+1},\\
\|L'(|W|^p|U|^q)\|_3 & \le C\|W\|_4^p\|U\|_3^q(T+R)^{p+q},\\
\|L'(|U|^r)\|_3 & \le C\|U\|_3^r(T+R)^{r+1}.
\end{array}
\right.
\end{equation}
\end{prop}
The proof of Proposition \ref{prop:apriori_zero}
is established in the next section.

\par
Let us proceed the proof of Theorem \ref{thm:lower-bound_zero}.
Set
\[
\begin{array}{ll}
N:=&\d2^{p+q-1}\sum_{\gamma=0}^1AE
(\|u_t^0\|_\infty^{p-\gamma}\|u_{tx}^0\|_\infty^\gamma\|u_0\|_\infty^q+
\|u_t^0\|_\infty^p\|u^0\|_\infty^{q-\gamma}\|u_x^0\|_\infty^\gamma)\\
&\d+\sum_{\gamma=0}^12^{r-\gamma}BE\|u^0\|_\infty^{r-\gamma}\|u_x^0\|^\gamma,
\end{array}
\]
where $E$ is the one in (\ref{apriori_linear}).
We note that
\[
\left\{
\begin{array}{lll}
\|u^0\|_{\infty}
&\le\|f\|_{L^\infty(\mathbf{R})} + \|g\|_{L^1(\mathbf{R})}&<\infty,\\
\|u_t^0\|_{\infty},\ \|u^0_x\|_{\infty}
&\le \|f'\|_{L^\infty(\mathbf{R})}+\|g\|_{L^\infty(\mathbf{R})}&<\infty,\\
\|u^0_{tx}\|_{\infty}
&\le \|f''\|_{L^\infty(\mathbf{R})}+\|g'\|_{L^\infty(\mathbf{R})}&<\infty.
\end{array}
\right.
\]
Assume that
\[
0<\e\le1.
\]
The four quantities, $\e_{2i}\ (i=1,2,3,4)$, are defined in the following;
\begin{equation}
\label{epsilon21}
\begin{array}{ll}
\e_{21}:= \min
&\left[\left\{2^{p+q}AC(5N)^{p+q}N^{-1}(2R)^{p+q}\right\}^{-1/[\min\{p+q,r\}(p+q-1)]},\right.\\
&\quad\{2^{p+q}AE\|u_t^0\|_\infty^p(5N)^qN^{-1}(2R)^q\}^{-1/[(q-1)\min\{p+q,r\}+p]},\\
&\quad\{2^{p+q}AE\|u^0\|_\infty^q(5N)^pN^{-1}(2R)^p\}^{-1/[(p-1)\min\{p+q,r\}+q]},\\
&\quad\left.\left\{2^rBC(5N)^rN^{-1}(2R)^{r+1}\right\}^{-1/[\min\{p+q,r\}(r-1)]}\right]
\end{array}
\end{equation}
and
\begin{equation}
\label{epsilon22}
\e_{22}:=\min\{\e_{221},\e_{222},\e_{223}\},
\end{equation}
where
\[
\begin{array}{l}
\e_{221}:=\min\left[\left\{2^{q+5}3^{p-1}pAC(5N)^{p+q-1}(2R)^{p+q}\right\}^{-1/[\min\{p+q,r\}(p+q-1)]},
\right.\\
\quad\{2^{q+4}3^{p-1}pAE\|u_t^0\|_\infty^{p-1}(5N)^qN(2R)^{q+1}\}^{-1/(p-1)},\\
\quad\{2^{q+4}3^{p-1}pAE\|u^0\|_\infty^q(5N)^{p-1}(2R)^p\}^{-1/[\min\{p+q,r\}(p-1)]},\\
\quad\left.\{2^{q+4}3^{p-1}pAE\|u_t^0\|_\infty^{p-1}\|u^0\|_\infty^q(2R)\}^{-1/(p+q-1)}\right]
\end{array}
\]
and
\[
\begin{array}{l}
\e_{222}:=\min\left[\{2^{p+4}3^qqAC(5N)^{p+q-1}(2R)^{p+q}\}^{-1/[\min\{p+q,r\}(p+q-1)]},\right.\\
\quad\{2^{p+3}3^qqAE\|u^0\|_\infty^{q-1}(5N)^p(2R)^{p+1}\}^{-1/[p\min\{p+q,r\}+q-1]},\\
\quad\{2^{p+4}3^qqAE\|u_t\|_\infty^p(5N)^{q-1}(2R)^q\}^{-1/[(q-1)\min\{p+q,r\}+p]},\\
\quad\left.\{2^{p+3}3^qqAE\|u_t^0\|^p\|u^0\|_\infty^{q-1}(2R)\}^{-1(p+q-1)}\right],\\
\e_{223}:=\min\left[\{2^33^rrBC(5N)^{r-1}(2R)^{r+1}\}^{-1/(r-1)},\right.\\
\quad\left.\quad\{2^33^rBE\|u^0\|_\infty^{r-1}(2R)\}^{-1/(r-1)}\right].
\end{array}
\]
Moreover, we set
\begin{equation}
\label{epsilon23}
\e_{23}:=\min\{\e_{21},\e_{231},\e_{232},\e_{233},\e_{234}\},
\end{equation}
where $\e_{21}$ is the one in (\ref{epsilon21}) and $\e_{23i}\ (i=1,2,3,4)$ are defined by
\[
\begin{array}{l}
\e_{231}:=\left[\{2^{p+q-1}pAC(5N)^{p+q}N^{-1}(2R)^{p+q}\}^{-1/[\min\{p+q,r\}(p+q-1)]},\right.\\
\quad\{2^{p+q-1}qAC(5N)^{p+q}N^{-1}(2R)^{p+q}\}^{-1/[\min\{p+q,r\}(p+q-1)]},\\
\quad\left.\{2^{r-1}qAC(5N)^rN^{-1}(2R)^{r+1}\}^{-1/[\min\{p+q,r\}(r-1)]}\right],
\end{array}
\]
\[
\begin{array}{l}
\e_{232}:=\min\big[\\
\quad\{2^{p+q-1}3^2pAE\|u_{tx}^0\|_\infty(5N)^{p+q-1}N^{-1}(2R)^{p+q-1}\}^{-1/[(p+q-1)\min\{p+q,r\}+1]},\\
\quad\{2^{p+q-1}3^2pAE\|u^0\|_\infty^q(5N)^pN^{-1}(2R)^p\}^{-1/(p\min\{p+q,r\}+q)},\\
\quad\{2^{p+q-1}3^2pAE\|u_{tx}^0\|_\infty\|u^0\|_\infty^q(5N)^{p-1}N^{-1}(2R)^{p-1}\}^{-1/[(p-1)\min\{p+q,r\}+q+1]},\\
\quad\{2^{p+q-1}3^2pAE\|u_t^0\|_\infty^{p-1}(5N)^{q+1}N^{-1}(2R)^{q+1}\}^{-1/[(q+1)\min\{p+q,r\}+p-1]},\\
\quad\{2^{p+q-1}3^2pAE\|u_t^0\|_\infty^{p-1}\|u_{tx}^0\|_\infty(5N)^qN^{-1}(2R)^q\}^{-1/(q\min\{p+q,r\}+p)},\\
\quad\left.\{2^{p+q-1}3^2pAE\|u_t^0\|_\infty^{p-1}\|u^0\|_\infty^q5(2R)\}^{-1/(\min\{p+q,r\}+p+q-1)}\right],
\end{array}
\]
\[
\begin{array}{l}
\e_{233}:=\min\big[\\
\quad\{2^{p+q-1}3^2qAE\|u_{tx}^0\|_\infty(5N)^{p+q-1}N^{-1}(2R)^{p+q-1}\}^{-1/[(p+q-1)\min\{p+q,r\}+1]},\\
\quad\{2^{p+q-1}3^2qAE\|u^0\|_\infty^{q-1}(5N)^{p+1}N^{-1}(2R)^{p+1}\}^{-1/[(p+1)\min\{p+q,r\}+q-1]},\\
\quad\{2^{p+q-1}3^2qAE\|u^0\|_\infty^{q-1}\|u_x^0\|_\infty(5N)^pN^{-1}(2R)^p\}^{-1/(p\min\{p+q,r\}+q)},\\
\quad\{2^{p+q-1}3^2qAE\|u_t^0\|_\infty^p(5N)^qN^{-1}(2R)^q\}^{-1/(q\min\{p+q,r\}+p)},\\
\quad\{2^{p+q-1}3^2qAE\|u_t^0\|_\infty^p\|u_x^0\|_\infty(5N)^{q-1}N^{-1}(2R)^{q-1}\}^{-1/[(q-1)\min\{p+q,r\}+p+1]},\\
\quad\left.\{2^{p+q-1}3^2qAE\|u_t^0\|_\infty^p\|u^0\|_\infty^{q-1}5(2R)\}^{-1/(\min\{p+q,r\}+p+q-1)}\right]
\end{array}
\]
and
\[
\begin{array}{l}
\e_{234}:=\min\big[\\
\quad\{2^{r-1}3rBE\|u_x^0\|_\infty(5N)^{r-1}N^{-1}(2R)^{r-1}\}^{-1/[(r-1)\min\{p+q,r\}+1]},\\
\quad\left.\{2^{r-1}3rBE\|u^0\|_\infty^{r-1}5(2R)\}^{-1/(\min\{p+q,r\}+r-1)}\right].
\end{array}
\]
Finally, we also set
\begin{equation}
\label{epsilon24}
\begin{array}{ll}
\e_{24}:= \min
&\left[\{2^{p+q+2}pAC(5N)^{p+q-1}(2R)^{p+q}\}^{-1/[\min\{p+q,r\}(p+q-1)]},\right.\\
&\quad\{2^{p+q+2}pAE\|u_t^0\|_\infty^{p-1}(5N)^q(2R)^{q+1}\}^{-1/(q\min\{p+q,r\})},\\
&\quad\{2^{p+q+2}pAE\|u^0\|_\infty^q(5N)^q(2R)^{p-1}\}^{-1/[(p-1)\min\{p+q,r\}]},\\
&\quad\left.\{2^{p+q+2}pAE\|u_t^0\|_\infty^{p-1}\|u^0\|_\infty^q(2R)\}^{-1}\right].
\end{array}
\end{equation}

\vskip10pt
\par\noindent
{\bf The convergence of the sequence $\v{\{(U_j,W_j)\}}$.}
\par
It follows from (\ref{U_j,W_j}), Proposition \ref{prop:apriori_linear} and \ref{prop:apriori_zero}
that
\[
\begin{array}{ll}
\|U_{j+1}\|_3
&\le A\|L(|W_j + \e u^0_t|^p|U_j+\e u^0|^q)\|_3+ B\|L(|U_j + \e u^0|^r)\|_3\\
&\le 2^{p+q}A\|L\{(|W_j|^p+|\e u^0_t|^p)(|U_j|^q+|\e u^0|^q)\}\|_3\\
&\quad+2^rB\|L(|U_j |^r+|\e u^0|^r)\|_3\\
&\le 2^{p+q}A\left\{\|L(|W_j|^p|U_j|^q)\|_3+\e^p\|L(|u_t^0|^p|U_j|^q)\|_3\right.\\
&\qquad\left.+\e^q\|L(|u^0|^q|W_j|^p)\|_3+\e^{p+q}\|L(|u^0_t|^p|u^0|^q)\|_3\right\}\\
&\quad+2^rB\left\{\|L(|U_j|^r)\|_3+\e^q\|L(|u^0|^r)\|_3 \right\}\\
&\le 2^{p+q}A\left\{C\|W_j\|_4^p\|U_j\|^q(T+R)^{p+q}
+E\e^p\|u_t^0\|_\infty^p\|U_j\|_3^q(T+R)^q\right.\\
&\quad+\left.E\e^q\|u^0\|_\infty^q\|W_j\|_4^p(T+R)^p+E\e^{p+q}\|u_t^0\|_\infty^p\|u^0\|_\infty^q\right\}\\
&\quad+2^r B\left\{C\|U_j\|_3^r(T+R)^{r+1} + E\e^r\|u^0\|_\infty^r \right\}.\\
&\le 2^{p+q}AC\|W_j\|_4^p\|U_j\|_3^q(T+R)^{p+q}\\
&\quad+2^{p+q}AE\|u_t^0\|_\infty^p\e^p\|U_j\|_3^q(T+R)^q\\
&\quad+2^{p+q}AE\|u^0\|_\infty^q\e^q\|W_j\|_4^p(T+R)^p\\
&\quad+2^rBC\|U_j\|_3^r(T+R)^{r+1} +N\e^{\min\{p+q,r\}}
\end{array}
\]
and similarly
\[
\begin{array}{ll}
\|W_{j+1}\|_4
&\le A\|L'(|W_j + \e u_t^0 |^p|U_j+\e u^0|^q )\|_4+ B\|L'(|U_j + \e u_0|^r)\|_4\\
&\le 2^{p+q}AC\|W_j\|_4^p\|U_j\|_3^q(T+R)^{p+q}\\
&\quad+2^{p+q}AE\|u_t^0\|_\infty^p\e^p\|U_j\|_3^q(T+R)^q\\
&\quad+2^{p+q}AE\|u^0\|_\infty^q\e^q\|W_j\|_4^p(T+R)^p\\
&\quad+2^rBC\|U_j\|_3^r(T+R)^{r+1} +N\e^{\min\{p+q,r\}}
\end{array}
\]
Hence the boundedness of $\{(U_j,W_j)\}$, i.e.
\begin{equation}
\label{bound_(U,W)}
\|U_j\|_3,\|W_j\|_4\le 5N\e^{\min\{p+q,r\}}\quad(j\in\N),
\end{equation}
follows from
\begin{equation}
\label{condi5}
\left\{
\begin{array}{ll}
2^{p+q}AC(5N\e^{\min\{p+q,r\}})^{p+q}(T+R)^{p+q}
&\le N\e^{\min\{p+q,r\}},\\
2^{p+q}AE\|u_t^0\|_\infty^p\e^p(5N\e^{\min\{p+q,r\}})^q(T+R)^q
&\le N\e^{\min\{p+q,r\}},\\
2^{p+q}AE\|u^0\|_\infty^q\e^q(5N\e^{\min\{p+q,r\}})^p(T+R)^p
&\le N\e^{\min\{p+q,r\}},\\
2^rBC(5N\e^{\min\{p+q,r\}})^r(T+R)^{r+1}
&\le N\e^{\min\{p+q,r\}}.
\end{array}
\right.
\end{equation}
Since (\ref{epsilon21}) yields that
\[
\begin{array}{ll}
R \le C_1\min & \left\{\e^{-\min\{p+q,r\}(p+q-1)/(p+q)},
\e^{-[(q-1)\min\{p+q,r\}+p]/q},\right.\\
& \left.\quad\e^{-[(p-1)\min\{p+q,r\}+q]/p},\e^{-\min\{p+q,r\}(r-1)/(r+1)}\right\}
\end{array}
\]
for $0<\e\le\e_{21}$,
where
\begin{equation}
\label{C_1}
\begin{array}{ll}
\d C_1:=\frac{1}{2} \min
&\left[\{2^{p+q}AC(5N)^{p+q}N^{-1}\}^{-1/(p+q)},\right.\\
&\quad\{2^{p+q}AE\|u_t^0\|_\infty^p(5N)^qN^{-1}\}^{-1/q},\\
&\quad\{2^{p+q}AE\|u^0\|_\infty^q(5N)^pN^{-1}\}^{-1/p},\\
&\left.\quad\{2^rBC(5N)^rN^{-1}\}^{-1/(r+1)}\right],
\end{array}
\end{equation}
we find that (\ref{condi5}) as well as (\ref{bound_(U,W)}) follows from
\begin{equation}
\label{lifespan_step1}
\begin{array}{ll}
T\le C_1\min
& \left\{\e^{-\min\{p+q,r\}(p+q-1)/(p+q)},\e^{-[(q-1)\min\{p+q,r\}+p]/q},\right.\\
& \quad\left.\e^{-[(p-1)\min\{p+q,r\}+q]/p},\e^{-\min\{p+q,r\}(r-1)/(r+1)}\right\}
\end{array}
\end{equation}
for $0<\e\le\e_{21}$.
\par
Let us write down (\ref{lifespan_step1}) in each cases.
\begin{itemize}
\item
When $p+q\le(r+1)/2(<r)$,\\
$\min\{p+q,r\}=p+q$ and
\[
\frac{p+q-1}{p+q}-\frac{r-1}{r+1}=\frac{2(p+q)-(r+1)}{(p+q)(r+1)}\le0
\]
imply that (\ref{lifespan_step1}) is equivalent to
\[
T\le C_1\e^{-(p+q-1)}
\]
because of
\[
\frac{(q-1)(p+q)+p}{q}=\frac{(p-1)(p+q)-q}{p}=p+q-1.
\]
\item
When $(r+1)/2\le p+q\le r$,\\
$\min\{p+q,r\}=p+q$ and
\[
\frac{p+q-1}{p+q}-\frac{r-1}{r+1}=\frac{2(p+q)-(r+1)}{(p+q)(r+1)}\ge0
\]
imply that (\ref{lifespan_step1}) is equivalent to
\[
T\le C_1\e^{-(p+q)(r-1)/(r+1)}.
\]
\item
When $r\le p+q$,\\
$\min\{p+q,r\}=r$ implies that (\ref{lifespan_step1}) is equivalent to
\[
T\le C_1\e^{-r(r-1)/(r+1)}
\]
because of
\[
\frac{r(q-1)+p}{q}-\frac{r(r-1)}{r+1}=\frac{r(p+q-r)+qr+p-r}{q(r+1)}>0
\]
and similarly
\[
\frac{r(p-1)+q}{p}-\frac{r(r-1)}{r+1}>0.
\]
\end{itemize}

\par
Next, assuming (\ref{lifespan_step1}),
one can estimate $(U_{j+1}-U_j)$ and $(W_{j+1}-W_j)$ as follows.
It follows from the inequality
\[
\begin{array}{l}
\left||W_j+\e u_t^0|^p|U_j+\e u^0|^q-|W_{j-1}+\e u_t^0|^p|U_j+\e u^0|^q\right|\\
\le\left||W_j+\e u_t^0|^p-|W_{j-1}+\e u_t^0|^p\right||U_j+\e u^0|^q\\
\quad+|W_{j-1}+\e u_t^0|^p\left||U_j+\e u^0|^q-|U_{j-1}+\e u^0|^q\right|\\
\le 2^q3^pp(|W_j|^{p-1}+|W_{j-1}|^{p-1}+\e^{p-1}|u_t^0|^{p-1})|W_j-W_{j-1}|(|U_j|^q+\e^q|u^0|^q)\\
\quad+2^p3^qq(|W_{j-1}|^p+\e^p|u_t^0|^p)(|U_j|^{q-1}+|U_{j-1}|^{q-1}+\e^{q-1}|u|^{q-1})|U_j-U_{j-1}|,
\end{array}
\]
that
\[
\begin{array}{l}
\|U_{j+1}-U_j\|_3\\
\le
2^q3^{p-1}pA\|L\{(|W_{j-1}|^{p-1} + |W_j|^{p-1} + \e^{p-1}|u^0_t|^{p-1})\times\\
\qquad\times|W_j-W_{j-1}|(|U_j|^q+\e^q|u^0|^q)\}\|_3\\
\quad+2^p3^{q-1}qA\|L\{(|W_{j-1}|^p+\e^p|u_t^0|^p)\times\\
\qquad\times(|U_{j-1}|^{q-1} + |U_j|^{q-1} + \e^{q-1}|u^0|^{q-1})|U_j-U_{j-1}|\}\|_3 \\
\quad+3^{r-1}rB\|L\{(|U_{j-1}|^{r-1} + |U_j|^{r-1} + \e^{r-1}|u^0|^{r-1})|U_j-U_{j-1}|\}\|_3.
\end{array}
\]
Hence Propositions \ref{prop:apriori_linear} and \ref{prop:apriori_zero} yield that
\[
\begin{array}{l}
\|U_{j+1}-U_j\|_3\\
\le
2^q3^{p-1}pA\|U_j\|_3^q\{C\|W_{j-1}\|_4^{p-1}(T+R)^{p+q}+C\|W_j\|_4^{p-1}(T+R)^{p+q}\\
\qquad+ E\e^{p-1}\|u^0_t\|_\infty^{p-1}(T+R)^{q+1}\}\|W_j-W_{j-1}\|_3\\
\quad+2^q3^{p-1}pAE\e^q\|u^0\|_\infty^q\{\|W_{j-1}\|_4^{p-1}(T+R)^p+\|W_j\|_4^{p-1}(T+R)^p\\
\qquad+E\e^{p-1}\|u^0_t\|_\infty^{p-1}(T+R)\}\|W_j-W_{j-1}\|_3\\
\quad+2^p3^{q-1}qA\|W_{j-1}\|_4^p\{C\|U_{j-1}\|_3^{q-1}(T+R)^{p+q}+C\|U_j\|_3^{q-1}(T+R)^{p+q}\\
\qquad+ E\e^{q-1}\|u^0\|_\infty^{q-1}(T+R)^{p+1}\}\|U_j-U_{j-1}\|_3\\
\quad+2^p3^{q-1}qAE\e^p\|u_t^0\|_\infty^p\{\|U_{j-1}\|_3^{q-1}(T+R)^q+\|U_j\|_3^{q-1}(T+R)^q\\
\qquad+E\e^{q-1}\|u^0\|_\infty^{q-1}(T+R)\}\|U_j-U_{j-1}\|_3 \\
\quad+3^{r-1}rB\{C\|U_{j-1}\|_3^{r-1}(T+R)^{r+1}+C\|U_j\|_3^{r-1}(T+R)^{r+1}\\
\qquad+E\e^{r-1}\|u^0\|_\infty^{r-1}(T+R)\}\|U_j-U_{j-1}\|_3.
\end{array}
\]
Therefore (\ref{condi5}) implies that
\[
\begin{array}{l}
\|U_{j+1}-U_j\|_3\\
\le
\{2^{q+1}3^{p-1}pAC(5N\e^{\min\{p+q,r\}})^{p+q-1}(T+R)^{p+q}\\
\qquad+2^q3^{p-1}pAE\|u_t^0\|_\infty^{p-1}
(5N\e^{\min\{p+q.r\}})^q\e^{p-1}(T+R)^{q+1}\}\|W_j-W_{j-1}\|_3\\
\quad+\{2^{q+1}3^{p-1}pAE\|u^0\|_\infty^q\e^q(5N\e^{\min\{p+q,r\}})^{p-1}(T+R)^p\\
\qquad+2^q3^{p-1}pAE\|u_t^0\|_\infty^{p-1}\|u^0\|_\infty^q\e^{p-1}(T+R)\}\|W_j-W_{j-1}\|_3\\
\quad+\{2^{p+1}3^{q-1}qAC(5N\e^{\min\{p+q,r\}})^{p+q-1}(T+R)^{p+q}\\
\qquad+2^p3^{q-1}qAE\|u^0\|_\infty^{q-1}(5N\e^{\min\{p+q,r\}})^p\e^{q-1}(T+R)^{p+1}\}\|U_j-U_{j-1}\|_3\\
\quad+\{2^{p+1}3^{q-1}qAE\|u_t^0\|_\infty^p\e^p(5N\e^{\min\{p+q,r\}})^{q-1}(T+R)^q\\
\qquad+2^p3^{q-1}qAE\|u_t^0\|_\infty^p\|u^0\|_\infty^{q-1}\e^{p+q-1}(T+R)\}\|U_j-U_{j-1}\|_3 \\
\quad+\{2\cdot3^{r-1}rBC(5N\e^{\min\{p+q,r\}})^{r-1}(T+R)^{r+1}\\
\qquad+3^{r-1}rBE\|u^0\|_\infty^{r-1}\e^{r-1}(T+R)\}\|U_j-U_{j-1}\|_3.
\end{array}
\]
We note that $\|W_{j+1}-W_j\|_4$ has the same upper bound
as $\|U_{j+1}-U_j\|_3$ in view of Propositions \ref{prop:apriori_linear} and \ref{prop:apriori_zero}.
Here we employ H\"older's inequality like the one in the proof of Theorem
\ref{thm:lower-bound_non-zero}.

\par
Therefore the convergence of $\{(U_j,W_j)\}$ follows from
\begin{equation}
\label{convergence_zero}
\begin{array}{l}
\|U_{j+1}-U_j\|_3+\|W_{j+1}-W_j\|_4\\
\d\le\frac{1}{2}\left(\|U_j-U_{j-1}\|_3+\|W_j-W_{j-1}\|_4\right)
\end{array}
\quad(j\ge2)
\end{equation}
provided (\ref{condi5}) and
\begin{equation}
\label{condi6}
\begin{array}{l}
\left.
\begin{array}{l}
2^{q+1}3^{p-1}pAC(5N\e^{\min\{p+q,r\}})^{p+q-1}(T+R)^{p+q},\\
2^q3^{p-1}pAE\|u_t^0\|_\infty^{p-1}(5N\e^{\min\{p+q,r\}})^q\e^{p-1}(T+R)^{q+1},\\
2^{q+1}3^{p-1}pAE\|u^0\|_\infty^q\e^q(5N\e^{\min\{p+q,r\}})^{p-1}(T+R)^p,\\
2^q3^{p-1}pAE\|u_t^0\|_\infty^{p-1}\|u^0\|_\infty^q\e^{p+q-1}(T+R)
\end{array}
\right\}
\le\d\frac{1}{16},
\\
\left.
\begin{array}{l}
2^{p+1}3^{q-1}qAC(5N\e^{\min\{p+q,r\}})^{p+q-1}(T+R)^{p+q},\\
2^p3^{q-1}qAE\|u^0\|_\infty^{q-1}(5N\e^{\min\{p+q,r\}})^p\e^{q-1}(T+R)^{p+1},\\
2^{p+1}3^{q-1}qAE\|u_t^0\|_\infty^p\e^p(5N\e^{\min\{p+q,r\}})^{q-1}(T+R)^q,\\
2^p3^{q-1}qAE\|u_t^0\|_\infty^p\|u^0\|_\infty^{q-1}\e^{p+q-1}(T+R),\\
2\cdot3^{r-1}rBC(5N\e^{\min\{p+q,r\}})^{r-1}(T+R)^{r+1},\\
3^{r-1}rBE\|u^0\|_\infty^{r-1}\e^{r-1}(T+R)
\end{array}
\right\}
\le\d\frac{1}{24}.
\end{array}
\end{equation}
are fulfilled.
Since (\ref{epsilon22}) yields that
\[
\begin{array}{ll}
R \le 
C_2\min
&\left\{\e^{-\min\{p+q,r\}(p+q-1)/(p+q)},\e^{-(q\min\{p+q,r\}+p-1)/(q+1)},\right.\\
&\quad\e^{-[(p-1)\min\{p+q,r\}+q]/p},\e^{-(p+q-1)},\\
&\quad\e^{-(p\min\{p+q,r\}+q-1)/(p+1)},\e^{-[(q-1)\min\{p+q,r\}+p]/q},\\
&\quad\left.\e^{-\min\{p+q,r\}(r-1)/(r+1)},\e^{-(r-1)}\right\},
\end{array}
\]
for $0<\e\le\e_{22}$, where
\[
\begin{array}{ll}
C_2:=\d\frac{1}{2}\min
&\left[
\{2^{q+5}3^{p-1}pAC(5N)^{p+q-1}\}^{-1/(p+q)},\right.\\
&\quad\{2^{q+4}3^{p-1}pAE\|u_t^0\|_\infty^{p-1}(5N)^q\}^{-1/(q+1)},\\
&\quad\{2^{q+5}3^{p-1}pAE\|u^0\|_\infty^q(5N)^{p-1}\}^{-1/p},\\
&\quad\{2^{q+4}3^{p-1}pAE\|u_t^0\|_\infty^{p-1}\|u^0\|_\infty^q\}^{-1},\\
&\quad\{2^{p+4}3^qqAC(5N)^{p+q-1}\}^{-1/(p+q)},\\
&\quad\{2^{p+3}3^qqAE\|u^0\|_\infty^{q-1}(5N)^p\}^{-1/(p+1)},\\
&\quad\{2^{p+4}3^qqAE\|u_t^0\|_\infty^p(5N)^{q-1}\}^{-1/q},\\
&\quad\{2^{p+3}3^qqAE\|u_t^0\|_\infty^p\|u^0\|_\infty^{q-1}\}^{-1},\\
&\quad\{2^33^rrBC(5N)^{r-1}\}^{-1/{r+1}},\\
&\quad\left.\{2^33^rrBE\|u^0\|_\infty^{r-1}\}^{-1}\right],
\end{array}
\]
we find that (\ref{condi6}) as well as (\ref{convergence_zero}) follows from
\begin{equation}
\label{lifespan_step2}
\begin{array}{ll}
T\le C_2\min
&\left\{\e^{-\min\{p+q,r\}(p+q-1)/(p+q)},\e^{-(q\min\{p+q,r\}+p-1)/(q+1)},\right.\\
&\quad\e^{-[(p-1)\min\{p+q,r\}+q]/p},\e^{-(p+q-1)},\\
&\quad\e^{-(p\min\{p+q,r\}+q-1)/(p+1)},\e^{-[(q-1)\min\{p+q,r\}+p]/q},\\
&\quad\left.\e^{-\min\{p+q,r\}(r-1)/(r+1)},\e^{-(r-1)}\right\}
\end{array}
\end{equation}
for $0<\e\le\e_{22}$.
\par
Let us write down (\ref{lifespan_step2}) in each cases.
\begin{itemize}
\item
When $p+q\le(r+1)/2(<r)$,\\
$\min\{p+q,r\}=p+q$ and the same reason as the argument after (\ref{lifespan_step1})
imply that (\ref{lifespan_step2}) is equivalent to
\[
T\le C_2\e^{-(p+q-1)}
\]
because of
\[
\frac{q(p+q)+p-1}{q+1}=\frac{p(p+q)+q-1}{p+1}=p+q-1.
\]
\item
When $(r+1)/2\le p+q\le r$,\\
similarly to the argument above, we find that (\ref{lifespan_step2}) is equivalent to
\[
T\le C_2\e^{-(p+q)(r-1)/(r+1)}.
\]
\item
When $r\le p+q$,\\
$\min\{p+q,r\}=r$ and the same reason as the argument after (\ref{lifespan_step1})
imply that (\ref{lifespan_step2}) is equivalent to
\[
T\le C_2\e^{-r(r-1)/(r+1)}
\]
because of
\[
\frac{qr+p-1}{q+1}-\frac{r(r-1)}{r+1}=\frac{r(p+q-r)+qr+p-1}{(q+1)(r+1)}>0
\]
and similarly
\[
\frac{pr+q-1}{p+1}-\frac{r(r-1)}{r+1}>0.
\]
\end{itemize}

\vskip10pt
\par\noindent
{\bf The convergence of the sequence $\v{\{((U_j)_x,(W_j)_x)\}}$.}
\par
Assume (\ref{lifespan_step1}) and (\ref{lifespan_step2}).
Then we have (\ref{bound_(U,W)}) and (\ref{convergence_zero}).
Since it follows from (\ref{U_j_x,W_j_x}) that
\[
\begin{array}{ll}
|(U_{j+1})_x|
&\le 2^{p+q}A\{L'(|W_j|^p|U_j|^q) +\e^pL'(|u_t|^0|^p|U_j|^q)\\
&\qquad+\e^qL'(|W_j|^p|u^0|^q)+ \e^{p+q}L'(|u^0_t|^p|u^0|^q) \}\\
&\quad +2^rB\{L'(|U_j|^r + \e^rL'(|u^0|^r)\}
\end{array}
\]
and
\[
\begin{array}{ll}
|(W_{j+1})_x|
&\le 2^{p+q-1}pA\{L'(|W_j|^{p-1}|(W_j)_x||U_j|^q)+ \e L'(|W_j|^{p-1}|u_{tx}^0||U_j|^q)\\
&\qquad+\e^qL'(|W_j|^{p-1}|(W_j)_x||u^0|^q)+ \e^{q+1}L'(|W_j|^{p-1}|u_{tx}^0||u^0|^q)\\
&\qquad+\e^{p-1}L'(|u_t^0|^{p-1}|(W_j)_x||U_j|^q)+\e^pL'(|u_t^0|^{p-1}|u_{tx}^0||U_j|^q)\\
&\qquad+\e^{p+q-1}L'(|u_t^0|^{p-1}|(W_j)_x||u^0|^q)+\e^{p+q}L'(|u_t^0|^{p-1}|u_{tx}^0||u^0|^q)\}\\
&\quad+2^{p+q-1}qA\{L'(|W_j|^p|U_j|^{q-1}|(U_j)_x|)+\e L'(|W_j|^p|U_j|^{q-1}|u_x^0|)\\
&\qquad+\e^{q-1}L'(|W_j|^p|u^0|^{q-1}|(U_j)_x|)+\e^qL'(|W_j|^p|u^0|^{q-1}|u_x^0|)\\
&\qquad+\e^pL'(|u_t^0|^p|U_j|^{q-1}|(U_j)_x|)+\e^{p+1}L'(|u_t^0|^p|U_j|^{q-1}|u_x^0|)\\
&\qquad+\e^{p+q-1}L'(|u_t^0|^p|u^0|^{q-1}|(U_j)_x|)+\e^{p+q}L'(|u_t^0|^p|u^0|^{q-1}|u_x^0|)\}\\
&\quad+2^{r-1}rB\{L'(|U_j|^{r-1}|(U_j)_x|) + \e L'(|U_j|^{r-1}|u_x^0|)\\
&\qquad+\e^{r-1}L'(|u^0|^{r-1}|(U_j)_x| + \e^rL'(|u^0|^{r-1}|u^0_x|))\},
\end{array}
\]
Proposition \ref{prop:apriori_linear} and \ref{prop:apriori_zero} yield that
\[
\begin{array}{ll}
\|(U_{j+1})_x\|_3
&\le 2^{p+q}A\{C\|W_j\|_4^p\|U_j\|_3^q(T+R)^{p+q}\\
&\qquad+\e^p E\|u_t^0\|_\infty^p\|U_j\|_3^q(T+R)^q\\
&\qquad+\e^qE\|W_j\|_4^p\|u^0\|_\infty^q(T+R)^p\\
&\qquad+ \e^{p+q}E\|u^0_t\|_\infty^p\|u^0\|_\infty^q)\}\\
&\quad +2^rB\{C\|U_j\|_3^r(T+R)^{r+1} + \e^rE\|u^0\|_\infty^r)\}
\end{array}
\]
and
\[
\|(W_{j+1})_x\|_4\le Z_1+Z_2+Z_3,
\]
where $Z_i\ (i=1,2,3)$ are defined by
\[
\begin{array}{ll}
Z_1:= & 2^{p+q-1}pA\{C\|W_j\|_4^{p-1}\|(W_j)_x\|_4\|U_j\|_3^q(T+R)^{p+q}\\
&+\e E\|u_{tx}^0\|_\infty\|W_j\|_4^{p-1}\|U_j\|_3^q(T+R)^{p+q-1}\\
&+\e^qE\|u^0\|_\infty^q\|W_j\|_4^{p-1}\|(W_j)_x\|_4(T+R)^p\\
&+ \e^{q+1}E\|u_{tx}^0\|_\infty\|u^0\|_\infty^q\|W_j\|_4^{p-1}(T+R)^{p-1}\\
&+\e^{p-1}E\|u_t^0\|_\infty^{p-1}\|(W_j)_x\|_4\|U_j\|_3^q(T+R)^{q+1}\\
&+\e^pE\|u_t^0\|_\infty^{p-1}\|u_{tx}^0\|_\infty\|U_j\|_3^q(T+R)^q\\
&+\e^{p+q-1}E\|u_t^0\|_\infty^{p-1}\|u^0\|_\infty^q\|(W_j)_x\|_4(T+R)\\
&+\e^{p+q}E\|u_t^0\|_\infty^{p-1}\|u_{tx}^0\|_\infty\|u^0\|_\infty^q\},
\end{array}
\]
\[
\begin{array}{ll}
Z_2:= & 2^{p+q-1}qA\{C\|W_j\|_4^p\|U_j\|_3^{q-1}\|(U_j)_x\|_3(T+R)^{p+q}\\
&+\e E\|u_x^0\|_\infty\|W_j\|_4^p\|U_j\|_3^{q-1}(T+R)^{p+q-1}\\
&+\e^{q-1}E\|u^0\|_\infty^{q-1}\|W_j\|_4^p\|(U_j)_x\|_3(T+R)^{p+1}\\
&+\e^qE\|u^0\|_\infty^{q-1}\|u_x^0\|_\infty\|W_j\|_4^p(T+R)^p\\
&+\e^pE\|u_t^0\|_\infty^p\|U_j\|_3^{q-1}\|(U_j)_x\|_3(T+R)^q\\
&+\e^{p+1}E\|u_t^0\|_\infty^p\|u_x^0\|_\infty\|U_j\|_3^{q-1}(T+R)^{q-1}\\
&+\e^{p+q-1}E\|u_t^0\|_\infty^p\|u^0\|_\infty^{q-1}\|(U_j)_x\|_3(T+R)\\
&+\e^{p+q}E\|u_t^0\|_\infty^p\|u^0\|_\infty^{q-1}\|u_x^0\|_\infty\},
\end{array}
\]
and
\[
\begin{array}{ll}
Z_3:=
&2^{r-1}rB\{C\|U_j\|_3^{r-1}\|(U_j)_x\|_3(T+R)^{r+1}\\
&+\e E\|u_x^0\|_\infty\|U_j\|_3^{r-1}(T+R)^{r-1}\\
&+\e^{r-1}E\|u^0\|_\infty^{r-1}\|(U_j)_x\|(T+R)\\
&\e^rE\|u^0\|_\infty^{r-1}\|u^0_x\|_\infty\}.
\end{array}
\]
Hence the boundedness of $\left\{\left((U_j)_x,(W_j)_x\right)\right\}$, i.e.
\begin{equation}
\label{bound_(U_x,W_x)}
\|(U_j)_x\|_3,\|(W_j)_x\|_4\le5N\e^{\min\{p+q,r\}}\quad (j\in\mathbf{N}),
\end{equation}
follows from (\ref{condi5}) for the $(U_j)_x$-component,
\begin{equation}
\label{condi7}
\left\{
\begin{array}{ll}
2^{p+q-1}pAC(5N\e^{\min\{p+q,r\}})^{p+q}(T+R)^{p+q} & \le N\e^{\min\{p+q,r\}},\\
2^{p+q-1}qAC(5N\e^{\min\{p+q,r\}})^{p+q}(T+R)^{p+q}\ & \le N\e^{\min\{p+q,r\}},\\
2^{r-1}rBC(5N\e^{\min\{p+q,r\}})^r(T+R)^{r+1} & \le N\e^{\min\{p+q,r\}}
\end{array}
\right.
\end{equation}
and the following three groups of conditions for the $(W_j)_x$-component;
\begin{equation}
\label{condi8-1}
\begin{array}{l}
\left.
\begin{array}{l}
\e E\|u_{tx}^0\|_\infty(5N\e^{\min\{p+q,r\}})^{p+q-1}(T+R)^{p+q-1},\\
\e^qE\|u^0\|_\infty^q(5N\e^{\min\{p+q,r\}})^p(T+R)^p,\\
\e^{q+1}E\|u_{tx}^0\|_\infty\|u^0\|_\infty^q(5N\e^{\min\{p+q,r\}})^{p-1}(T+R)^{p-1},\\
\e^{p-1}E\|u_t^0\|_\infty^{p-1}(5N\e^{\min\{p+q,r\}})^{q+1}(T+R)^{q+1},\\
\e^pE\|u_t^0\|_\infty^{p-1}\|u_{tx}^0\|_\infty(5N\e^{\min\{p+q,r\}})^q(T+R)^q,\\
\e^{p+q-1}E\|u_t^0\|_\infty^{p-1}\|u^0\|_\infty^q(5N\e^{\min\{p+q,r\}})(T+R)
\end{array}
\right\}\\
\d\le\frac{2}{3}\cdot\frac{1}{6}\cdot\frac{1}{2^{p+q-1}pA}N\e^{\min\{p+q,r\}}
,
\end{array}
\end{equation}
\begin{equation}
\label{condi8-2}
\begin{array}{l}
\left.
\begin{array}{l}
\e E\|u_x^0\|_\infty(5N\e^{\min\{p+q,r\}})^{p+q-1}(T+R)^{p+q-1},\\
\e^{q-1}E\|u^0\|_\infty^{q-1}(5N\e^{\min\{p+q,r\}})^{p+1}(T+R)^{p+1},\\
\e^qE\|u^0\|_\infty^{q-1}\|u_x^0\|_\infty(5N\e^{\min\{p+q,r\}})^p(T+R)^p,\\
\e^pE\|u_t^0\|_\infty^p(5N\e^{\min\{p+q,r\}})^q(T+R)^q,\\
\e^{p+1}E\|u_t^0\|_\infty^p\|u_x^0\|_\infty(5N\e^{\min\{p+q,r\}})^{q-1}(T+R)^{q-1},\\
\e^{p+q-1}E\|u_t^0\|_\infty^p\|u^0\|_\infty^{q-1}(5N\e^{\min\{p+q,r\}})(T+R)
\end{array}
\right\}\\
\d\le\frac{2}{3}\cdot\frac{1}{6}\cdot\frac{1}{2^{p+q-1}qA}N\e^{\min\{p+q,r\}}

\end{array}
\end{equation}
and
\begin{equation}
\label{condi8-3}
\begin{array}{l}
\left.
\begin{array}{l}
\e E\|u_x^0\|_\infty(5N\e^{\min\{p+q,r\}})^{r-1}(T+R)^{r-1},\\
\e^{r-1}E\|u^0\|_\infty^{r-1}(5N\e^{\min\{p+q,r\}})(T+R)
\end{array}
\right\}
\\
\d\le\frac{2}{3}\cdot\frac{1}{2}\cdot\frac{1}{2^{r-1}rB}N\e^{\min\{p+q,r\}},
\end{array}
\end{equation}

\par
Note that (\ref{epsilon23}) yields
\[
\begin{array}{ll}
R\le C_3\min
&
\left\{\e^{-\min\{p+q,r\}(p+q-1)/(p+q)},\e^{-\min\{p+q,r\}(r-1)/(r+1)},\right.\\
&\quad\e^{-[(p+q-2)\min\{p+q,r\}+1]/(p+q-1)},\e^{-[(p-1)\min\{p+q,r\}+q]/p},\\
&\quad\e^{-[(p-2)\min\{p+q,r\}+q+1]/(p-1)},\e^{-[q\min\{p+q,r\}+p-1]/(q+1)},\\
&\quad\e^{-[(q-1)\min\{p+q,r\}+p]/q},\e^{-(p+q-1)},\\
&\quad\e^{-[p\min\{p+q,r\}+q-1]/(p+1)},\e^{-[(q-2)\min\{p+q,r\}+p+1]/(q-1)}\\
&\quad\left.\e^{-[(r-2)\min\{p+q,r\}+1]/(r-1)},\e^{-(r-1)}\right\}
\end{array}
\]
for $0<\e\le\e_{23}$. Here we set
\[
C_3:=\frac{1}{2}\min\{2C_1,C_{31},C_{32},C_{33},C_{34}\},
\]
where $C_1$ is the one in (\ref{C_1}) and $C_{3i}\ (i=1,2,3,4)$ are defined by
\[
\begin{array}{ll}
C_{31}:=\min
&\left[\{2^{p+q-1}pAC(5N)^{p+q}N^{-1}\}^{-1/(p+q)},\right.\\
&\quad\{2^{p+q-1}qAC(5N)^{p+q}N^{-1}\}^{-1/(p+q)},\\
&\quad\left.\{2^{r-1}rBC(5N)^rN^{-1}\}^{-1/(r+1)}\right],
\end{array}
\]
\[
\begin{array}{ll}
C_{32}:=\min
&\left[\{2^{p+q-1}3^2pAE\|u_{tx}^0\|_\infty(5N)^{p+q-1}N^{-1}\}^{-1/(p+q-1)},\right.\\
&\quad\{2^{p+q-1}3^2pAE\|u^0\|_\infty^q(5N)^pN^{-1}\}^{-1/p},\\
&\quad\{2^{p+q-1}3^2pAE\|u_{tx}^0\|_\infty\|u^0\|_\infty^q(5N)^{p-1}N^{-1}\}^{-1/(p-1)},\\
&\quad\{2^{p+q-1}3^2pAE\|u_t^0\|_\infty^{p-1}(5N)^{q+1}N^{-1}\}^{-1/(q+1)},\\
&\quad\left.\{2^{p+q-1}3^2pAE\|u_t^0\|_\infty^{p-1}\|u^0\|_\infty^q5\}^{-1}\right],
\end{array}
\]
\[
\begin{array}{ll}
C_{33}:=\min
&\left[\{2^{p+q-1}3^2qAE\|u_x^0\|_\infty(5N)^{p+q-1}N^{-1}\}^{-1/(p+q-1)},\right.\\
&\quad\{2^{p+q-1}3^2qAE\|u^0\|_\infty^{q-1}(5N)^{p+1}N^{-1}\}^{-1/(p+1)},\\
&\quad\{2^{p+q-1}3^2qAE\|u^0\|_\infty^{q-1}\|u_x^0\|_\infty(5N)^pN^{-1}\}^{-1/p},\\
&\quad\{2^{p+q-1}3^2qAE\|u_t^0\|_\infty^p(5N)^qN^{-1}\}^{-1/q},\\
&\quad\{2^{p+q-1}3^2qAE\|u_t^0\|_\infty^p\|u_x^0\|_\infty(5N)^{q-1}N^{-1}\}^{-1/(q-1)},\\
&\quad\left.\{2^{p+q-1}3^2qAE\|u_t^0\|_\infty^p\|u^0\|_\infty^{q-1}5\}^{-1}\right]
\end{array}
\]
and
\[
\begin{array}{ll}
C_{34}:=\min
&\left[\{2^{r-1}3rBE\|u_x^0\|_\infty(5N)^{r-1}N^{-1}\}^{-1/(r-1)},\right.\\
&\quad\left.\{2^{r-1}3rBE\|u^0\|_\infty^{r-1}5\}^{-1}\right].
\end{array}
\]
Therefore we find that (\ref{condi7}), (\ref{condi8-1}), (\ref{condi8-2}) and (\ref{condi8-3})
as well as (\ref{bound_(U_x,W_x)}) follow from
\begin{equation}
\label{lifespan_step3}
\begin{array}{ll}
T\le C_3&\min
\left\{\e^{-\min\{p+q,r\}(p+q-1)/(p+q)},\e^{-\min\{p+q,r\}(r-1)/(r+1)},\right.\\
&\quad\e^{-[(p+q-2)\min\{p+q,r\}+1]/(p+q-1)},\e^{-[(p-1)\min\{p+q,r\}+q]/p},\\
&\quad\e^{-[(p-2)\min\{p+q,r\}+q+1]/(p-1)},\e^{-[q\min\{p+q,r\}+p-1]/(q+1)},\\
&\quad\e^{-[(q-1)\min\{p+q,r\}+p]/q},\e^{-(p+q-1)},\\
&\quad\e^{-[p\min\{p+q,r\}+q-1]/(p+1)},\e^{-[(q-2)\min\{p+q,r\}+p+1]/(q-1)}\\
&\quad\left.\e^{-[(r-2)\min\{p+q,r\}+1]/(r-1)},\e^{-(r-1)}\right\}
\end{array}
\end{equation}
for $0<\e\le\e_{23}$.

\par
Let us write down this inequality in each case.
First we note that
\[
\begin{array}{l}
\d\frac{(\gamma-1)\min\{p+q,r\}+p+q-\gamma}{\gamma}\\
\d=\min\{p+q,r\}-1+\frac{p+q-\min\{p+q,r\}}{\gamma}\\
\d\ge\min\{p+q,r\}-1+\frac{p+q-\min\{p+q,r\}}{p+q}\\
\d=\frac{(p+q-1)\min\{p+q,r\}}{p+q}
\end{array}
\]
holds for any $\gamma$ satisfying $0<\gamma\le p+q$.
Therefore (\ref{lifespan_step3}) can be diminished as
\[
\begin{array}{ll}
T\le C_3\min
&\left\{\e^{-\min\{p+q,r\}(p+q-1)/(p+q)},\e^{-\min\{p+q,r\}(r-1)/(r+1)},\right.\\
&\quad\left.\e^{-[(r-2)\min\{p+q,r\}+1]/(r-1)},\e^{-(r-1)}\right\}
\end{array}
\]
\begin{itemize}
\item
When $(2<)p+q\le(r+1)/2(<r)$,\\
$\min\{p+q,r\}=p+q$ and the same reason as the argument after (\ref{lifespan_step1})
imply that (\ref{lifespan_step3}) is equivalent to
\[
T\le C_3\e^{-(p+q-1)}
\]
because of
\[
\begin{array}{l}
\d\frac{(r-2)(p+q)+1}{r-1}-\frac{(p+q)(r-1)}{r+1}\\
\d=\frac{(r^2-r-2)(p+q)+r+1-(p+q)(r^2-2r+1)}{(r-1)(r+1)}\\
\d=\frac{(r-3)(p+q)+r+1}{(r-1)(r+1)}\ge\frac{(r-2)(p+q)+1}{(r-1)(r+1)}>0.
\end{array}
\]
\item
When $(r +1)/2\le p+q \le r$,
similarly to the argument above, we find that (\ref{lifespan_step3}) is equivalent to
\[
T\le C_3\e^{-(p+q)(r-1)/(r+1)}.
\]
\item
When $p+q\ge r$,
$\min\{p+q,r\}=r$ and the argument after (\ref{lifespan_step1})
imply that (\ref{lifespan_step3}) is equivalent to
\[
T\le C_3\e^{-r(r-1)/(r+1)}
\]
because of $\min\{p+q,r\}=r$ and
\[
\begin{array}{l}
\d\frac{(r-2)r+1}{r-1}-\frac{r(r-1)}{r+1}\\
\d=\frac{(r-3)r+r+1}{(r-1)(r+1)}=\frac{r-1}{r+1}>0.
\end{array}
\]
\end{itemize}
\par
Next, assuming (\ref{lifespan_step1}), (\ref{lifespan_step2}) and (\ref{lifespan_step3}),
we shall estimate $\{(U_{j+1})_x-(U_j)_x\}$ and $\{(W_{j+1})_x-(W_j)_x\}$.
First we deal with $\{(U_{j+1})_x-(U_j)_x\}$.
It is easy to see
\[
\begin{array}{l}
|(U_{j+1})_x - (U_j)_x|\\
\le L'(A||W_j + \e u^0_t|^p|U_j+\e u^0|^q - |W_{j-1} + \e u^0_t|^p||U_{j-1}+\e u^0|^q|)\\
\quad+L'(B||U_j + \e u^0|^r - |U_{j-1} + \e u^0|^r|)\\
\le 2^q3^{p-1}pAL'\{(|W_{j-1}|^{p-1} + |W_j|^{p-1} + \e^{p-1}|u^0_t|^{p-1})\times\\
\qquad\times(|U_j|^q+\e^q|u^0|^q)|W_j-W_{j-1}|\}\\
\quad+2^p3^{q-1}qAL'\{(|W_j|^p+\e^p|u_t^0|^p)(|U_j|^{q-1}+|U_{j-1}|^{q-1}+\e^{q-1}|u^0|^{q-1})\times\\
\qquad\times|U_j-U_{j-1}|\}\\
\quad+3^{r-1}rBL'\{(|U_{j-1}|^{r-1} + |U_j|^{r-1} + \e^{q-1}|u^0|^{r-1})|U_j - U_{j-1}|\}.
\end{array}
\]
Hence it follows from Propositions \ref{prop:apriori_linear} and  \ref{prop:apriori_zero} that
\[
\begin{array}{l}
\|(U_{j+1})_x - (U_j)_x\|_3\\
\le 2^q3^{p-1}pAC(\|W_{j-1}\|_4^{p-1} + \|W_j\|_4^{p-1})\|U_j\|_3^q\|W_j - W_{j-1}\|_4(T+R)^{p+q}\\
\quad+2^q3^{p-1}pAE\e^{p-1}\|u_t^0\|_\infty^{p-1}\|U_j\|_3^q\|W_j-W_{j-1}\|_4(T+R)^{q+1}\\
\quad+2^q3^{p-1}pA\e^q\|u^0\|_\infty^q(\|W_{j-1}\|_4^{p-1}+\|W_j\|_4^{p-1})\|W_j-W_{j-1}\|_4(T+R)^p\\
\quad+2^q3^{p-1}pA\e^{p+q-1}\|u_t^0\|_\infty^{p-1}\|u^0\|_\infty^q\|W_j-W_{j-1}\|_4(T+R)\\
\quad+2^p3^{q-1}qAC\|W_j\|_4^p(\|U_j\|_3^{q-1}+\|U_{j-1}\|_3^{q-1})\|U_j-U_{j-1}\|_3(T+R)^{p+q}\\
\quad+2^p3^{q-1}qAE\e^{q-1}\|u^0\|_\infty^{q-1}\|W_j\|_4^p\|U_j-U_{j-1}\|_3(T+R)^{p+1}\\
\quad+2^p3^{q-1}qAE\e^p\|u_t^0\|_\infty^p(\|U_j\|_3^{q-1}+\|U_{j-1}\|_3^{q-1})\|U_j-U_{j-1}\|_3(T+R)^q\\
\quad+2^p3^{q-1}qAE\e^{p+q-1}\|u^0_t\|_\infty^p\|u^0\|_\infty^{q-1}\|W_j - W_{j-1}\|_4(T+R)\\
\quad+3^{r-1}rBC(\|U_{j-1}\|_3^{r-1} + \|U_j\|_3^{r-1} )\|U_j - U_{j-1}\|_3(T+R)^{r+1}\\
\quad+3^{r-1}rBE\e^{r-1}\|u^0\|_\infty^{r-1}\|U_j - U_{j-1}\|_3(T+R)
\end{array}
\]
Since (\ref{lifespan_step1}) and (\ref{lifespan_step2}) yield (\ref{convergence_zero}),
we have that
\[
\|U_{j+1}-U_j\|_3 + \|W_{j+1} - W_{j}\|_3\le O\left(\frac{1}{2^j}\right)
\quad\mbox{as}\ j\rightarrow\infty.
\]
This fact implies that
\begin{equation}
\label{convergence_U_x}
\|(U_{j+1})_x-(U_j)_x\|_3 = O\left(\frac{1}{2^j}\right)
\quad\mbox{as}\ j\rightarrow\infty.
\end{equation}

\par
Next, we estimate $\{(W_{j+1})_x-(W_j)_x\}$.
It is also easy to see
\[
\begin{array}{l}
|(W_{j+1})_x - (W_j)_x|\\
\le pAL'\{||W_j + \e u^0_t|^{p-2}(W_j + \e u^0_t)((W_j)_x + \e u^0_{tx})|U_j+\e u^0|^q\\
\qquad-|W_{j-1} + \e u^0_t|^{p-2}(W_{j-1} + \e u^0_t)((W_{j-1})_x + \e u^0_{tx})|U_{j-1}+\e u^0|^q|\}\\
\quad+qAL'\{||W_j+\e u_t^0|^p|U_j + \e u^0|^{q-2}(U_j + \e u^0)((U_j)_x + \e u^0_x)\\
\qquad-|W_{j-1}+\e u_t^0|^p|U_{j-1} + \e u^0|^{q-2}(U_{j-1} + \e u^0)((U_{j-1})_x + \e u^0_x)|\}\\
\quad+rBL'\{||U_j + \e u^0|^{r-2}(U_j + \e u^0)((U_j)_x + \e u^0_x)\\
\qquad-|U_{j-1} + \e u^0|^{r-2}(U_{j-1} + \e u^0)((U_{j-1})_x + \e u^0_x)|\}.
\end{array}
\]
In order to estimate the quantities in the right hand side of this inequality, we employ
\[
\begin{array}{l}
||W_j + \e u^0_t|^{p-2}(W_j + \e u^0_t)((W_j)_x + \e u^0_{tx})|U_j+\e u^0|^q\\
-|W_{j-1} + \e u^0_t|^{p-2}(W_{j-1} + \e u^0_t)((W_{j-1})_x + \e u^0_{tx})|U_{j-1}+\e u^0|^q|\\
\le||W_j + \e u^0_t|^{p-2}(W_j + \e u^0_t)-|W_{j-1} + \e u^0_t|^{p-2}(W_{j-1} + \e u^0_t)|\times\\
\quad \times|(W_j)_x + \e u^0_{tx}||U_j+\e u^0|^q\\
\quad+|W_{j-1} + \e u^0_t|^{p-1}|(W_j)_x - (W_{j-1})_x||U_j+\e u^0|^q\\
\quad+|W_{j-1}+\e u_t^0|^{p-1}|(W_{j-1})_x+\e u_{tx}|||U_j+\e u^0|^q-|U_{j-1}+\e u^0|^q|
\end{array}
\]
and
\[
\begin{array}{l}
||W_j + \e u^0_t|^{p-2}(W_j + \e u^0_t)- |W_{j-1} + \e u^0_t|^{p-2}(W_{j-1} + \e u^0_t)|\\
\le\left\{
\begin{array}{ll}
\begin{array}{l}
3^{p-2}(p-1)(|W_j|^{p-2}+|W_{j-1}|^{p-2} + \e^{p-2}|u^0_t|^{p-2})\times\\
\quad\times|W_j - W_{j-1}|
\end{array}
& \mbox{for}\ p\ge 2,\\
2|W_j - W_{j-1}|^{p-1} & \mbox{for}\ 1< p < 2.
\end{array}
\right.
\end{array}
\]
Hence the same manner as  estimating $\{(U_{j+1})_x-(U_j)_x\}$ yields that
\[
\begin{array}{l}
\|(W_j)_x - (W_{j-1})_x\|_4\\
\le 2^{p+q-1}pA\|L'\{(|W_{j-1}|^{p-1}+\e^{p-1}|u_t^0|^{p-1})\times\\
\quad\times(|U_j|^q+\e^q|u^0|^q)|(W_j)_x-(W_{j-1})|\}\|_4
+\d O\left(\frac{1}{2^{j\min\{p-1,q-1,r-1,1\}}}\right)
\end{array}
\]
as $j\rightarrow\infty$.
Then it follows from (\ref{lifespan_step1}) as well as (\ref{bound_(U,W)}) that
\[
\begin{array}{ll}
2^{p+q-1}pA\|L'\{(|W_{j-1}|^{p-1}+\e^{p-1}|u^0_t|^{p-1})\times\\
\quad\times(|U_j|^q+\e^q|u^0|^q)|(W_j)_x - (W_{j-1})_x|)\|_4
\\
\le 2^{p+q-1}pAC\|W_{j-1}\|_4^{p-1}\|U_j\|_3^q\|(W_j)_x-(W_{j-1})_x\|_4(T+R)^{p+q}\\
\quad+2^{p+q-1}pAE\e^{p-1}\|u_t^0\|_\infty^{p-1}\|U_j\|_3^q\|(W_j)_x-(W_{j-1})_x\|_4(T+R)^{q+1}\\
\quad+2^{p+q-1}pAE\e^q\|u^0\|_\infty^q\|W_{j-1}\|_4^{p-1}\|(W_j)_x-(W_{j-1})_x\|_4(T+R)^p\\
\quad+2^{p+q-1}pAE\e^{p+q-1}\|u_t^0\|_\infty^{p-1}\|u^0\|_\infty^q
\|(W_j)_x-(W_{j-1})_x\|_4(T+R).
\end{array}
\]
Therefore we obtain that
\begin{equation}
\label{convergence_W_x}
\|(W_{j+1})_x-(W_j)_x\|_4
\le \frac{1}{2}\|(W_j)_x-(W_{j-1})_x\|_4
+O\left(\frac{1}{2^{j\min\{p-1,q-1,1\}}}\right)
\end{equation}
as $j\rightarrow\infty$, which shows the convergence of $\{(W_j)_x\}$,
provided
\begin{equation}
\label{condi9}
\left.
\begin{array}{ll}
2^{p+q-1}pAC(5N\e^{\min\{p+q,r\}})^{p+q-1}(T+R)^{p+q},\\
2^{p+q-1}pAE\e^{p-1}\|u_t^0\|_\infty^{p-1}(5N\e^{\min\{p+q,r\}})^q(T+R)^{q+1},\\
2^{p+q-1}pAE\e^q\|u^0\|_\infty^q(5N\e^{\min\{p+q,r\}})^{p-1}(T+R)^p,\\
2^{p+q-1}pAE\e^{p+q-1}\|u_t^0\|_\infty^{p-1}\|u^0\|_\infty^q(T+R)
\end{array}
\right\}
\le\frac{1}{8}
\end{equation}
holds.
Since (\ref{epsilon24}) yields that
\[
\begin{array}{ll}
R\le C_4\min & \{\e^{-\min\{p+q,r\}(p+q-1)/(p+q)},\e^{-[q\min\{p+q,r\}+p-1]/(q+1)},\\
& \quad\e^{-[(p-1)\min\{p+q,r\}+q]/p},\e^{-(p+q-1)}\}
\end{array}
\]
for $0<\e\le \e_{24}$, where
\[
\begin{array}{ll}
C_4:=\d\frac{1}{2} \min & \left[\{2^{p+q+2}pAC(5N)^{p+q-1}\}^{-1/(p+q)},\right.\\
&\quad\{2^{p+q+2}pAE\|u_t^0\|_\infty^{p-1}(5N)^q\}^{-1/(q+1)},\\
&\quad\{2^{p+q+2}pAE\|u\|_\infty^q(5N)^{p-1}\}^{-1/p},\\
&\quad\left.\{2^{p+q+2}pAE\|u_t^0\|_\infty^{p-1}\|u^0\|_\infty\}^{-1}\right]
\end{array}
\]
we find that (\ref{condi9}) as well as (\ref{convergence_U_x}) and (\ref{convergence_W_x})
follows from
\begin{equation}
\label{lifespan_step4}
\begin{array}{ll}
T\le C_4\min & \{\e^{-\min\{p+q,r\}(p+q-1)/(p+q)},\e^{-[q\min\{p+q,r\}+p-1]/(q+1)},\\
& \quad\e^{-[(p-1)\min\{p+q,r\}+q]/p},\e^{-(p+q-1)}\}
\end{array}
\end{equation}
for $0<\e\le \e_{24}$.
We note that (\ref{lifespan_step4}) is a part of (\ref{lifespan_step2})
for which $C_4$ is replace with $C_2$.

\vskip10pt
\par\noindent
{\bf Continuation of the proof.}
\par
The convergence of the sequence $\{(U_j,W_j)\}$ to $(U,W)$ in the closed subspace of $Y$
satisfying
\[
 \|U\|_3,\|(U_x)\|_3,\|W\|_4,\|(W)_x\|_4\le 5N\e^{\min\{p+q,r\}}
 \]
is established by  (\ref{epsilon21}), (\ref{epsilon22}), (\ref{epsilon23}), (\ref{epsilon24}),
(\ref{lifespan_step1}), (\ref{lifespan_step2}), (\ref{lifespan_step3}) and (\ref{lifespan_step4}).
Therefore the statement of Theorem \ref{thm:lower-bound_zero} is established with
\[
\left\{
\begin{array}{l}
\d c=\min\left\{C_1,C_2,C_3,C_4\right\},\\
\e_2=\min\{1,\e_{21},\e_{22},\e_{23},\e_{24}\}.
\end{array}
\right.
\]
\hfill$\Box$


\section{Proof of Proposition \ref{prop:apriori_zero}}

\par
In this section, we prove a priori estimate (\ref{apriori_zero}).
Note that three estimates with $|U|^r$ are already obtained by Proposition 5.2
in Morisawa, Sasaki and Takamura \cite{MST},
so that we shall prove other three estimates with $|W|^p|U|^q$.
Here a positive constant $C$ independent of $T$ and $\e$
may change from line to line.

\par
It follows from the assumption on the supports and the definition of $L$ that
\[
|L(|W|^p|U|^q)(x,t)|\le C\|W\|_4^p\|U\|_3^qJ(x,t)
\quad\mbox{for}\ |x|\le t+R,
\]
where we set
\[
\begin{array}{ll}
J(x,t):=&\d\int_0^tds\int_{x-t+s}^{x+t-s}
\{\chi_D(s,y)+(1-\chi_D(s,y))(s+|y|+R)^p\}\times\\
&\quad\times(s+|y|+R)^q\chi_{\mbox{\footnotesize supp}(U,W)}(s,y)dy.
\end{array}
\]
First, we consider the case of $x\ge0$.
From now on, we employ the change of variables
\begin{equation}
\label{characteristic}
\alpha=s+y,\ \beta=s-y.
\end{equation}
For $(x,t)\in D$, extending the domain of the integral, we have that
\[
\begin{array}{ll}
J(x,t)
&\d\le C\int_{-R}^Rd\beta\int_{-R}^{x+t}(\alpha+R)^{p+q}d\alpha\\
&\d\quad+C\int_R^{t-x}(\beta+R)^{p+q}d\beta\int_{-R}^Rd\alpha\\
&\d\quad+C\int_R^{t-x}d\beta\int_R^{t+x}(\alpha+R)^qd\alpha\\
&\le C(t+x+R)^{p+q+1}\\
&\quad+C(t-x+R)^{p+q+1}\\
&\quad+C(t+x+R)^{q+1}(t-x+R)\\
&\le C(T+R)^{p+q}(t+x+R).
\end{array}
\]
For $t+x\ge R$ and $|t-x|\le R$, we also have that
\[
\begin{array}{ll}
J(x,t)
& \d\le C+C\int_{-R}^{t-x}d\beta\int_R^{t+x}(\alpha+R)^{p+q}d\alpha\\
& \d\le C(T+R)^{p+q}(t+x+R).
\end{array}
\]
For $t+x\le R$, it is trivial that
\[
J(x,t)\le C.
\]
Summing up, we obtain that
\[
\begin{array}{c}
|L(|W|^p|U|^q)(x,t)|\le C\|W\|_4^p\|U\|_3^q(T+R)^{p+q}(t+x+R)\\
\mbox{for}\ 0\le x\le t+R.
\end{array}
\]
The case of $x\le0$ is similar to the one above, so we omit the details.
Therefore we obtain the first inequality in (\ref{apriori_zero}).

\par
Next, we shall show the third inequality in (\ref{apriori_zero}).
It follows from the assumption on the supports and the definition of $L'$ that
\[
\begin{array}{l}
|L'(|W|^p|U|^q)(x,t)|\\
\quad\le C\|W\|_4^p\|U\|_3^q\{J_+(x,t)+J_-(x,t)\}
\end{array}
\mbox{for}\ |x|\le t+R,
\]
where the integrals $J_+$ and $J_-$ are defined by
\[
\begin{array}{ll}
J_\pm(x,t):=&
\d\int_0^t
\{\chi_{I\pm}(x,t;s)+(1-\chi_{I\pm}(x,t;s))(s+|t-s\pm x|+R)^p\}\times\\
&\quad\times\chi_{E\pm}(x,t;s)(s+|t-s\pm x|+R)^qds
\end{array}
\]
and the characteristic functions $\chi_{I+},\chi_{I-},\chi_{E+}$ and $\chi_{E-}$ are defined by
\[
\begin{array}{l}
\chi_{I\pm}(x,t;s):=\chi_{\{s: s-|t-s\pm x|\ge R\}},\\
\chi_{E\pm}(x,t;s):=\chi_{\{s: |t-s\pm x|\le s+R\}}
\end{array}
\]
respectively.
First we note that it is sufficient to  estimate $J_\pm$ for $x\ge0$ due to its symmetry,
\[
J_+(-x,t)=J_-(x,t).
\]
For $(x,t)\in D\cap\{x\ge0\}$, we have
\[
\begin{array}{ll}
J_+(x,t)
&\d\le C\int_{(t+x-R)/2}^{(t+x+R)/2}(t+x+R)^{p+q}ds\\
&\d\quad+C\int_{(t+x+R)/2}^t(t+x+R)^qds\\
&\le C(t+x+R)^{p+q}+C(t+x+R)^{q+1}\\
&\le C(T+R)^{p+q}
\end{array}
\]
and
\[
\begin{array}{ll}
J_-(x,t)
&\d\le C\int_{(t-x-R)/2}^{(t-x+R)/2}(s+|t-s-x|+R)^{p+q}ds\\
&\d\quad+C\int_{(t-x+R)/2}^t(s+|t-s-x|+R)^qds\\
&\le C(t+x+R)^{p+q}+C(t+x+R)^{q+1}\\
&\le C(T+R)^{p+q}.
\end{array}
\]
For $t+x\ge R$ and $|t-x|\le R$, we have
\[
J_+(x,t)
\le C\int_{(t+x-R)/2}^t(t+x+R)^{p+q}ds
\le C(T+R)^{p+q}(t+x+R)\]
and
\[
\begin{array}{ll}
J_-(x,t)
&\d\le C\int_0^t(s+|t-s-x|+R)^{p+q}ds\\
&\d\le C\int_0^{t-x}(t-x+R)^{p+q}ds
+C\int_{t-x}^t(2s-t+x+R)^{p+q}ds\\
&\le C(t-x+R)^{p+q+1}+C(t+x+R)^{p+q+1}\\
&\le C(T+R)^{p+q}(t+x+R).
\end{array}
\]
It is trivial that $J_{\pm}(x,t)\le C$ for $t+x\le R$.
Therefore we obtain the third inequality in (\ref{apriori_zero}).

\par
The fifth inequality in (\ref{apriori_zero})
readily follows from the computations above.
The proof of Proposition \ref{prop:apriori_zero} is now completed.
\hfill$\Box$


\section{Proofs of Theorem \ref{thm:upper-bound_non-zero}
and Theorem \ref{thm:upper-bound_zero}}

\par
The essential argument to obtain the upper bound of the lifespan is the following.
Let $u$ be a classical solution of (\ref{IVP_gcombined}) in a time interval $[0,T]$.
If $T$ is bigger than some quantity depending on $\e$,
we will meet a contradiction to the fact that
$u$ is a classical solution.
This situation gives us that the lifespan should be less than the quantity
due to its definition.

\vskip10pt
\par\noindent
{\bf Proof of  Theorem \ref{thm:upper-bound_non-zero}.}
\par
Neglecting the second term of our equation as
\[
A|u_t|^p|u|^q+B|u|^r\ge A|u_t|^p|u|^q,
\]
in (\ref{IVP_gcombined}), we have that there exists a constant
$\e_{31}=\e_{31}(f,g,p,q,A,R)$ such that the contradiction appears provided
\begin{equation}
\label{lifespan_zhou1}
T>C_{31}\e^{-(p+q-1)}\quad\mbox{if}\ \int_{\R}g(x)dx>0
\end{equation}
holds for $0<\e\le\e_{31}$ and some positive constant $C_{31}$ independent of $\e$.
Because it is already obtained by Zhou \cite{Zhou} for the equation
\[
u_{tt}-u_{xx}=|u_t|^p|u|^q.
\]
in which it is trivial that \lq\lq$=$" can be replaced with \lq\lq$\ge A\times$".
As stated in Introduction, we shall repeat its proof in Appendix below.

\par
By virtue of the same reason and Zhou \cite{Zhou92}, making use of
\[
A|u|_t^p|u|^q+B|u|^r\ge B|u|^r,
\]
 we have that there exists a constant
$\e_{32}=\e_{32}(f,g,r,B,R)$ such that the contradiction appears provided
\begin{equation}
\label{lifespan_zhou2}
T>C_{32}\e^{-(r-1)/2}\quad\mbox{if}\ \int_{\R}g(x)dx>0
\end{equation}
holds for $0<\e\le\e_{32}$ and some positive constant $C_{32}$ independent of $\e$.
Therefore, taking
\[
\e_3=\min\{\e_{31},\e_{32},1\}\quad\mbox{and}\quad C=\min\{C_{31},C_{32}\},
\]
we have the desired lifespan estimate by (\ref{lifespan_zhou1}) and (\ref{lifespan_zhou2}).
\hfill$\Box$

\vskip10pt
\par\noindent
{\bf Proof of  Theorem \ref{thm:upper-bound_zero}.}
\par
First we shall prove that there exists a constant
$\e_{41}=\e_{41}(f,p,q,A,R)$ such that the contradiction appears provided
\begin{equation}
\label{lifespan_new1}
T>C_{41}\e^{-(p+q-1)}\quad\mbox{if (\ref{positive_zero}) is fulfilled}
\end{equation}
holds for $0<\e\le\e_{41}$ and some positive constant $C_{41}$ independent of $\e$.
To this end, define the blow-up set
\begin{equation}
\label{Sigma}
\Sigma:=\left\{(x,t)\in\R\times[0,T]\ :\ x\ge0,\ t+x\ge R,\ 0<t-x<\frac{R}{2}\right\}.
\end{equation}
Then, it follows from the positiveness of the nonlinear term and
representation of $u,u_t$ in (\ref{u}), (\ref{u_t}) that
\[
u(x,t)\ge\frac{\e}{2}\{f(x+t)+f(x-t)\}
\]
and
\[
u_t(x,t)\ge\frac{\e}{2}\{f'(x+t)-f'(x-t)\}.
\]
Hence the assumption on the initial data, (\ref{positive_zero}), implies that
\begin{equation}
\label{1st}
u(x,t),u_t(x,t)\ge\frac{1}{2}f_0\e\quad\mbox{in}\ \Sigma.
\end{equation}
\par
From now on, we employ the routine iteration procedure.
Assume an estimate
\begin{equation}
\label{estimate}
\left\{
\begin{array}{ll}
u(x,t) & \ge M_n(t+x-R)^{a_n}(t-x)^{b_n},\\
u_t(x,t)\ & \ge M_n(t+x-R)^{a_n}(t-x)^{c_n}
\end{array}
\right.
\mbox{for}\ (x,t)\in\Sigma
\end{equation}
holds, where $a_n,b_n,c_n\ge0$ and $M_n>0$.
All the sequences $\{a_n\},\{b_n\},\ \{c_n\}$ and $\{M_n\}$ are defined later.
Then it follows from (\ref{u}), (\ref{u_t}) and (\ref{characteristic}) that
\[
\begin{array}{ll}
u(x,t)
&\d\ge\frac{AM_n^{p+q}}{4}\int_0^{t-x}\beta^{qb_n+pc_n}d\beta
\int_R^{t+x}(\alpha-R)^{(p+q)a_n}d\alpha\\
&\d\ge\frac{AM_n^{p+q}}{4\{(p+q)a_n+1\}\{qb_n+pc_n+1\}}\times\\
&\quad\times(t+x-R)^{(p+q)a_n+1}(t-x)^{qb_n+pc_n+1}
\end{array}
\]
and
\[
\begin{array}{ll}
u_t(x,t)
&\d\ge\frac{A}{2}\int_{(t-x+R)/2}^t|u_t(x-t+s,s)|^p|u(x-t+s,s)|^qds\\
&\d\ge\frac{AM_n^{p+q}}{2}(t-x)^{qb_n+pc_n}\int_{(t-x+R)/2}^t(2s-t+x-R)^{(p+q)a_n}ds\\
&\d\ge\frac{AM_n^{p+q}}{4\{(p+q)a_n+1\}}(t+x-R)^{(p+q)a_n+1}(t-x)^{qb_n+pc_n}
\end{array}
\]
for $(x,t)\in\Sigma$.
Hence (\ref{estimate}) holds for all $n\in\N$ provided
\[
\left\{
\begin{array}{ll}
a_{n+1}=(p+q)a_n+1,& a_1=0,\\
b_{n+1}=qb_n+pc_n+1, & b_1=0,\\
c_{n+1}=qb_n+pc_n, & c_1=0
\end{array}
\right.
\]
and
\[
M_{n+1}\le\frac{AM_n^{p+q}}{4\{(p+q)a_n+1\}\{(p+q)(b_n+c_n)+1\}},\ M_1=\frac{1}{2}f_0\e.
\]
It is easy to see that
\[
a_n=b_n+c_n=\frac{(p+q)^{n-1}-1}{p+q-1}\quad(n\in\N),
\]
which implies
\[
\begin{array}{l}
\{(p+q)a_n+1\}\{(p+q)(b_n+c_n)+1\}\\
\d\le\{(p+q)a_n+1\}^2=a_{n+1}^2\le\frac{(p+q)^{2n}}{(p+q-1)^2}.
\end{array}
\]
Therefore $M_n$ should be defined by
\[
M_{n+1}=AC_5(p+q)^{-2n}M_n^{p+q},\ M_1=\frac{1}{2}f_0\e,
\]
where we set
\[
C_5:=\frac{(p+q-1)^2}{4}>0,
\]
so that (\ref{estimate}) implies that
\begin{equation}
\label{lower-bound}
\begin{array}{l}
u(x,t)u_t(x,t)\\
\ge M_n^2(t+x-R)^{2a_n}(t-x)^{a_n}\\
\ge C_6\{(t+x-R)^2(t-x)\}^{-1/(p+q-1)}\exp\left\{Z(x,t)(p+q)^{n-1}\right\}
\end{array}
\end{equation}
for $(x,t)\in\Sigma$, where
\[
\begin{array}{ll}
Z(x,t)
&\d:=\frac{1}{p+q-1}\log\{(t+x-R)^2(t-x)\}\\
&\d\quad+\frac{2}{p+q-1}\log(AC_5)-4S_{p+q}\log(p+q)+2\log\left(\frac{1}{2}f_0\e\right),\\
C_6
&\d:=\exp\left(-\frac{2}{p+q-1}\log (AC_5)\right)>0.
\end{array}
\]
Indeed, $M_n$ satisfies
\[
\log M_{n+1}=\log(AC_5)-2n\log(p+q)+(p+q)\log M_n,
\]
which implies
\[
\begin{array}{l}
\log M_{n+1}\\
=\{1+(p+q)+\cdots+(p+q)^{n-1}\}\log(AC_5)\\
\quad-2\{n+(p+q)(n-1)+\cdots+(p+q)^{n-1}(n-(n-1))\}\log(p+q)\\
\quad+(p+q)^n\log M_1\\
\d=\frac{(p+q)^n-1}{p+q-1}\log(AC_5)-2(p+q)^{n-1}\log(p+q)\sum_{j=0}^{n-1}\frac{j+1}{(p+q)^j}\\
\quad+(p+q)^n\log M_1\\
\d\ge-\frac{1}{p+q-1}\log(AC_5)\\
\d\quad+(p+q)^n\left\{\frac{1}{p+q-1}\log(AC_5)-2S_{p+q}\log(p+q)+\log M_1\right\},
\end{array}
\]
where we set
\[
S_r:=\sum_{j=0}^\infty\frac{j+1}{r^{j+1}}<\infty.
\]

\par
In view of (\ref{lower-bound}), if there exists a point $(x_0,t_0)\in\Sigma$ such that
$Z(x_0,t_0)>0$, we have a contradiction $u(x_0,t_0)u_t(x_0,t_0)=\infty$ to the fact that
$u$ is a classical solution on the time interval $[0,T]$ with $T\ge t_0$
of (\ref{IVP_gcombined}) by letting $n\rightarrow\infty$.
Let us set
\[
t_0=x_0+\frac{R}{4}\quad\mbox{and}\quad t_0\ge\frac{5R}{4}(>1).
\]
Then, since we have
\[
(t_0+x_0-R)^2(t_0-x_0)\ge\frac{R}{4}t_0^2,
\]
$Z(x_0,t_0)>0$ follows from
\[
t_0^2>\frac{4}{R}\cdot\frac{(p+q)^{4(p+q-1)S_{p+q}}}{(AC_5)^2(2^{-1}f_0)^{2(p+q-1)}}\e^{-2(p+q-1)}.
\]
Therefore this inequality completes the proof of (\ref{lifespan_new1})
with
\[
C_{41}=\frac{2}{\sqrt{R}}\cdot\frac{(p+q)^{2(p+q-1)S_{p+q}}}{AC_5(2^{-1}f_0)^{p+q-1}}
\]
and $\e_{41}>0$ satisfying
\[
C_{41}\e_{41}^{-(p+q-1)}=\frac{5R}{4}.
\]


\par
Now we note that there exists a constant $\e_{42}=\e_{42}(f,r,B,R)$ such that
a contradiction appears provided
\begin{equation}
\label{lifespan_new2}
T>C_{42}\e^{-r(r-1)/(r+1)}\quad\mbox{if (\ref{positive_zero}) is fulfilled}
\end{equation}
holds for $0<\e\le\e_{42}$ and some positive constant $C_{42}$ independent of $\e$.
The assumption (\ref{positive_zero}) is stronger than the one of
Theorem 5.1 in Takamura \cite{Takamura15}
for the equation $u_{tt}-u_{xx}=|u|^r$.
Its proof is available also for $u_{tt}-u_{xx}\ge B|u|^r$ with a trivial modification.


\par
Finally, we shall deal with the case of the generalized combined effect, namely
\[
\frac{r+1}{2}\le p+q\le r.
\]
Set
\[
F(t):=\int_{\R}u(x,t)dx.
\]
Then the equation in (\ref{IVP_gcombined}) and the support of the solution (\ref{support_sol})
yield that
\[
F''(t)=\int_{\R}\{A|u_t(x,t)|^p|u(x,t)|^q+B|u(x,t)|^r\}dx
\]
Neglecting the first term in the integrand and making use of H\"older's inequality,
we have that
\begin{equation}
\label{ODI}
F''(t)\ge B\int_{\R}|u(x,t)|^rdx\ge 2^{1-r}B(t+R)^{-(r-1)}|F(t)|^r
\quad\mbox{for}\ t\ge0.
\end{equation}
On the other hand, neglecting the second term of the integrand, we have that
\[
F''(t)\ge A\int_{\R}|u_t(x,t)|^p|u(x,t)|^qdx
\ge A\int_{t-R/2}^t|u_t(x,t)|^p|u(x,t)|^qdx
\]
for $t\ge R/2$.
Then it follows from (\ref{1st}) that
\[
F''(t)\ge \frac{ARf_0^{p+q}}{2^{p+q+1}}\e^{p+q}
\quad\mbox{for}\ t\ge\frac{R}{2}.
\]
Integrating this inequality and employing the fact that
$F''(t)\ge0$ for $t\ge0$ and 
\[
\left\{
\begin{array}{ll}
F'(t) & \d\ge F'(0)=\e\int_{\R}g(x)dx=0,\\
F(t) & \d\ge F(0)=\e\int_{\R}f(x)dx>0
\end{array}
\right.
\mbox{for}\ t\ge0,
\]
we have that
\[
F'(t)\ge\frac{ARf_0^{p+q}}{2^{p+q+2}}\e^{p+q}t
\quad\mbox{for}\ t\ge R
\]
which yields that
\begin{equation}
\label{initial}
F(t)\ge\frac{ARf_0^{p+q}}{2^{p+q+4}}\e^{p+q}t^2
\quad\mbox{for}\ t\ge 2R.
\end{equation}

\par
We are now in a position to employ Lemma 2.2 in Takamura \cite{Takamura15}.
The assumption in (2.9) in \cite{Takamura15} is satisfied by (\ref{initial}) and setting
\[
\frac{ARf_0^{p+q}}{2^{p+q+4}}\e^{p+q}t_0^2=2F(0)=2\e\int_{\R}f(x)dx.
\]
Also with
\[
p=r,\ q=r-1,\ a=2,
\]
the blow-up condition (2.1) in \cite{Takamura15} implies the trivial inequality
\[
M=\frac{r-1}{2}\cdot2-\frac{r-1}{2}+1=\frac{r+1}{2}>0.
\]
If we put
\[
T_0=C_0'\left(\frac{ARf_0^{p+q}}{2^{p+q+4}}\e^{p+q}\right)^{-(r-1)/(r+1)}
\]
where $C_0'=C_0'(r,B)>0$ is the $C_0$ in Lemma 2.2 in \cite{Takamura15},
we have
\[
T_0\ge\max\{t_0,2R\}
\]
provided
\[
C_0'\left(\frac{ARf_0^{p+q}}{2^{p+q+4}}\e^{p+q}\right)^{-(r-1)/(r+1)}
\ge
\left(\frac{2^{p+q+5}}{ARf_0^{p+q}}\int_{\R}f(x)dx\e^{-(p+q-1)}\right)^{1/2}\ge 2R
\]
holds.
This inequality can be possible to be established for small $\e$ by the fact that
\[
\frac{p+q-1}{2}<(p+q)\frac{r-1}{r+1}
\]
is equivalent to
\[
(p+q)\frac{3-r}{r+1}<1.
\]
This inequality is trivial when $r\ge3$ and also follows from a trivial inequality
\[
r\frac{3-r}{r+1}<1
\]
for $1<r<3$ in this case, $(r+1)/2\le p+q\le r$.
Therefore it is possible to take $T_2=T_0$ in Lemma 2.2 in \cite{Takamura15},
so that there exists a constant $\e_{43}=\e_{42}(f,p,q,r,A,B,R)$ such that
the contradiction appears provided
\begin{equation}
\label{lifespan_new3}
T>C_{43}\e^{-r(r-1)/(r+1)}\quad\mbox{if (\ref{positive_zero}) is fulfilled}
\end{equation}
holds for $0<\e\le\e_{43}$ and some positive constant $C_{43}$ independent of $\e$.
According to the computations above, we can choose
\[
\e_{43}=\min\{\e_{431},\e_{432}\},
\]
where $\e_{431}$ and $\e_{432}$ are defined by
\[
C_0'\left(\frac{ARf_0^{p+q}}{2^{p+q+4}}\e_{431}^{p+q}\right)^{-(r-1)/(r+1)}
=
\left(\frac{2^{p+q+5}}{ARf_0^{p+q}}\int_{\R}f(x)dx\e_{431}^{-(p+q-1)}\right)^{1/2}
\]
and
\[
\left(\frac{2^{p+q+5}}{ARf_0^{p+q}}\int_{\R}f(x)dx\e_{432}^{-(p+q-1)}\right)^{1/2}=2R.
\]
Also it is possible to set
\[
C_{43}=2^{4/(r+1)}C_0'\left(\frac{ARf_0^{p+q}}{2^{p+q+4}}\right)^{-(r-1)/(r+1)}>0.
\]

\par
Recall Remark \ref{rem:relation}.
Summing up (\ref{lifespan_new1}), (\ref{lifespan_new2}) and (\ref{lifespan_new3}),
we obtain the statement of Theorem \ref{thm:upper-bound_zero}
by taking
\[
\e_4=\min\{1,\e_{41},\e_{42},\e_{43}\}
\quad\mbox{and}\quad
C=\min\{C_{41},C_{42},C_{43}\}.
\]
The proof of Theorem \ref{thm:upper-bound_zero} is completed now. 
\hfill$\Box$


\section*{Appendix}
As stated in Introduction, we repeat here the result and its proof,
restricted in one space dimension although,
of the unpublished paper by Zhou \cite{Zhou}, for the sake of the completeness of this paper.

\begin{thm}[Zhou \cite{Zhou}]
Let $A>0$ and $B=0$.
Assume (\ref{supp_initial}) and (\ref{positive_non-zero}).
Then, there exists a positive constant $\e_5=\e_5(f,g,p,q,A,R)>0$ such that
the lifespan $T(\e)$ of a classical solution of (\ref{IVP_gcombined})
satisfies
\begin{equation}
\label{lifespan_zhou}
T(\e)\le C\e^{-(p+q-1)}
\end{equation}
where $0<\e\le\e_5$, and $C$ is a positive constant independent of $\e$.
\end{thm}
\par\noindent
{\bf Proof.} The proof is almost same as Zhou \cite{Zhou01}
in which only the nonlinear term $|u_t|^p$ is considered.

\par
Let $u$ be a classical solution of (\ref{IVP_gcombined})
in the time interval $[0,T]$.
Set
\[
G:=\frac{1}{2}\int_{\R}g(x)dx>0.
\]
From now on, we restrict ourselves in the domain
\[
D':=\{(x,t)\in\R\times[0,T]\ :\ x\ge R\quad\mbox{and}\quad t-x\ge R\}.
\]
Then, it follows from (\ref{supp_initial})  and (\ref{u}) that
\[
u= G\e+AL(|u_t|^p|u|^q)\quad\mbox{in}\ D'.
\]
Inverting the order of variables and diminishing the domain in $L$, we have
\[
L(v)(x,t)
\ge\frac{1}{2}\int_R^xdy\int_{y-R}^{y+R}v(y,s)ds
\quad\mbox{for}\ (x,t)\in D'
\]
for any non-negative function $v=v(x,t)$
Applying $v=|u_t|^p|u|^q$ to the equation above and making use of (\ref{product}),
we obtain that
\begin{equation}
\label{fundamental}
P(x)\ge G\e+C_7\int_R^x|P(y)|^{p+q}dy\quad\mbox{for}\ x\ge R,
\end{equation}
where we set
\[
P(x):=u(x,x+R)\quad\mbox{and}\quad C_7:=\frac{A}{2}\left(\frac{p}{p+q}\right)^p(2R)^{1-p}>0.
\]
Because we have employed H\"older's inequality to have
\[
\int_{y-R}^{y+R}\left(|u(y,s)|^{(p+q)/p}\right)_sds
\le
\left\{\int_{y-R}^{y+R}\left|\left(|u(y,s)|^{(p+q)/p}\right)_s\right|^pds\right\}^{1/p}(2R)^{1-1/p}
\]
and (\ref{support_sol}) implies that
\[
\int_{y-R}^{y+R}\left(|u(y,s)|^{(p+q)/p}\right)_sds=|u(y,y+R)|^{(p+q)/p}.
\]

\par
Once (\ref{fundamental}) is obtained,
it is easy to reach to the desired conclusion by completely the same argument in \cite{Zhou01}
in which $p$ is replaced with $p+q$.
The proof is now completed.
\hfill$\Box$.

\section*{Acknowledgement}
\par
The second author is partially supported
by the Grant-in-Aid for Young Scientists (No.18K13447), 
Japan Society for the Promotion of Science.
The fourth author is partially supported
by the Grant-in-Aid for Scientific Research (A) (No.22H00097), 
Japan Society for the Promotion of Science.


\bibliographystyle{plain}

\end{document}